\documentclass[a4paper]{article}

\usepackage{amsmath,graphicx}
\usepackage{amsfonts}
\usepackage[english]{babel}
\usepackage[utf8x]{inputenc}
\usepackage{amsmath}
\usepackage{graphicx}
\usepackage[colorinlistoftodos]{todonotes}
\usepackage{subcaption}
\usepackage[bottom]{footmisc}
\usepackage{amsthm}
\usepackage{setspace}
\usepackage{adjustbox}
\usepackage{multirow}
\usepackage{hhline}
\usepackage{tablefootnote}
\usepackage{authblk}
\title{Higher-Order Block Term Decomposition \\
for Spatially Folded fMRI Data\thanks{The research leading to these results has received funding from the European Union's Seventh Framework Programme (H2020-MSCA-ITN-2014) under grant agreement No.~642685 MacSeNet.}}

\author[1,2]{Christos Chatzichristos}
\author[3,2]{Eleftherios Kofidis}
\author[4,2]{Yiannis Kopsinis}
\author[1,2]{Sergios Theodoridis}

\affil[1]{Department of Informatics and Telecommunications, University of Athens, Athens, Greece (stheodor@di.uoa.gr)}
\affil[2]{Computer Technology Institute \& Press ``Diophantus" (CTI), Greece (chatzichris@cti.gr)}
\affil[3]{Department of Statistics and Insurance Science, University of Piraeus, Piraeus, Greece (kofidis@unipi.gr)}
\affil[4]{LIBRA MLI Ltd, Edinburgh, UK (kopsinis@ieee.org)}

\date{}

\begin{document}

\maketitle

\begin{abstract}
The growing use of neuroimaging technologies generates a massive amount of biomedical data that exhibit high dimensionality. Tensor-based analysis of brain imaging data has been proved quite effective in exploiting their multiway nature. The advantages of tensorial methods over matrix-based approaches have also been demonstrated in the characterization of functional magnetic resonance imaging (fMRI) data, where the spatial (voxel) dimensions are commonly grouped (unfolded) as a single way/mode of the 3-rd order array, the other two ways corresponding to time and subjects. However, such methods are known to be ineffective in more demanding scenarios, such as the ones with strong noise and/or significant overlapping of activated regions. This paper aims at investigating the possible gains from a better exploitation of the spatial dimension, through a higher- (4 or 5) order tensor modeling of the fMRI signal. In this context, and in order to increase the degrees of freedom of the modeling process, a higher-order  Block Term Decomposition (BTD) is applied, for the first time in fMRI analysis. Its effectiveness is demonstrated via extensive simulation results. 
\end{abstract}

\section{Introduction}
\label{sec:intro}
Functional Magnetic Resonance Imaging (fMRI) is a noninvasive technique for studying brain activity, which receives an increasing attention in the last decade or so. During an fMRI experiment, a series of brain images is acquired, while the subject possibly performs a set of tasks responding to external stimuli. Changes in the measured blood-oxygen-level dependent (BOLD) signal are used to examine different types of activation in the brain. There are several objectives in the analysis of fMRI data, the most common of which are the localization of regions of the brain, that are activated by a task, and the determination of the functional brain connectivity~\cite{2008_lindquist_statistical,2004_huettel_functional}. 

The localization of the activated areas in the human brain is a challenging ``cocktail party" problem, where several people are talking (areas activated) simultaneously behind a wall (skull). Our goal is to distinguish those areas (spatial maps) as well as activation patterns (time courses) through some blind source separation (decomposition) method \cite{2015_theodoridis_machine,2006_calhoun_unmixing}. Each source is the outcome of a combination of a time course with a spatial map.

\begin{figure} [b]
\centering
\captionsetup{justification=centering}
\includegraphics[width=0.99\textwidth]{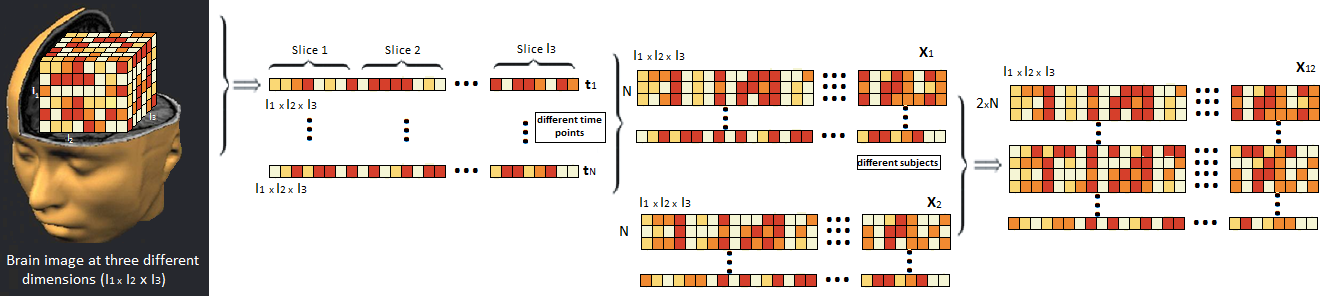}
\caption{Brain images are unfolded in vectors, which are then stacked in matrices.} \label{fig:concat}
\end{figure}

In fMRI studies of the brain function, the structure of the data involves multiple modes, such as trial, task condition, subject, in addition to the intrinsic dimensions of time and space \cite{2004_andersen_structure-seeking}. The use of multivariate bi-linear (matrix-based) methods through concatenating different dimensions, have been, so far, the state of the art for extracting information concerning spatial and temporal features in data from a single fMRI run. Such methods have also been extended to multi-subject experiments~\cite{2006_calhoun_unmixing,2012_andersen_identifying,2011_dea_iva,2011_penny_statistical,2008_xiong_ica}.

After acquiring an fMRI image at a time instance $n$ (Fig.~\ref{fig:concat}), the three-dimensional data (referred to here as folded data) is reshaped (unfolded) into a vector $\mathbf{t}_n$. $N$ of those vectors (fMRI images of the same subject, at different time instants) are stacked together in a matrix, $\mathbf{X}_1$. Similarly, the data coming from a second subject (multi-subject case) are concatenated into another matrix, $\mathbf{X}_2$,  and the two matrices are joined together to form $\mathbf{X}_{12}$; the latter will then be decomposed by some blind source separation method. Thus, the five-dimensional problem of a multi-subject fMRI experiment has been transformed into a two-dimensional problem. This reshaping (unfolding) of higher-order data (data of higher than two dimensions) into two-way arrays leads to decompositions that are inherently non-unique (due to rotation invariance), and, most importantly, can result in a loss of the multi-way (neighborhood) linkages and useful hidden (latent) components during the modeling process.

The multi-dimensional nature of the data (and any additional information found in it) can be retained in multilinear models, which, in general, produce representations that are unique~\cite{2000_sidiropoulos_uniqueness}, can improve the ability of extracting spatiotemporal modes of interest, and can facilitate subsequent interpretations that are neurophysiologically meaningful \cite{2004_andersen_structure-seeking}. In the example of Fig.~\ref{fig:concat}, instead of being concatenated into a bigger matrix $\mathbf{X}_{12}$, the matrices $\mathbf{X}_{1}$ and $\mathbf{X}_{2}$ can be used to form a third-order tensor (multi-way array) with the
dimensions being voxels $\times$ times $\times$ subjects, and hence a tensor decomposition method can be mobilized for the blind source separation task \cite{2015_cichocki_tensor}. Although the results from such tensor-based methods have been very promising, providing, in most of the cases, better spatial and temporal localization of the activity compared to the matrix-based approaches, the step of the initial spatial unfolding of the data to a vector $\mathbf{t}_n$ seems to have been ``inherited" from matrix-based approaches \cite{2004_andersen_structure-seeking,2005_beckmann_tensorial,2008_morup_shift-invariant,2013_davidson_network}, being more reminiscent of bi-linear rather than of fully multilinear modeling.

It is therefore of interest to explore the application of models of order higher than three and investigate possible benefits from fully exploiting the underlying spatial information. Indeed, it is demonstrated in this paper that, if properly exploited, such a higher-dimensional viewpoint can improve the localization of the sources and decrease the effect of noise. The measured fMRI signal is corrupted by random noise and various nuisance components, that arise due to both hardware reasons (thermal noise, scanner drift, etc.) as well as to the physiology-related nature of the experiments (subject's motion during the acquisition, heart beat, breathing, etc.). Furthermore, the noise is not homogeneous throughout the brain, which renders its suppression even more difficult \cite{2006_lund_non-white,2015_bright_is,2006_wink_bold}.

In order to use all the spatial information in our data, our kick off point will be to bypass the initial step of unfolding the data into $\mathbf{t}_n$ vectors. Furthermore, the Block-Term Decomposition (BTD) model~\cite{2008_de_lathauwer_decompositions-2,2008_de_lathauwer_decompositions-1,2008_de_lathauwer_decompositions,2011_de_lathauwer_blind,2012_de_lathauwer_block} will be adopted, for a first time in fMRI, in view of its higher modeling potential and its increased robustness to rank estimation errors. 

The goal of this paper is to propose a method that can fully exploit the original geometry of the problem in order to improve the accuracy and interpretability of the decomposition results. Through extensive simulation results, in scenarios well known in the literature, it will be demonstrated that the proposed method can overcome drawbacks of the state-of-the art techniques, improving the accuracy and interpretability of the decomposition results even in challenging scenarios. 

The rest of this paper is organized as follows. Section 2 briefly reviews the state-of-the art in tensorial fMRI analysis, pointing out the known limitations of the existing methods. The BTD-based method along with a uniqueness argument for the 4-way BTD are presented in Section 3. Section 4 reports and discusses a number of simulation results in multi-subject scenarios while in Section 5 the new possibilities offered by the proposed approach for tensor-based analysis of single-subject cases are pointed out. Conclusions are drawn in Section 6. 

\subsection{Notation}

Vectors, matrices and higher-order tensors are denoted by bold lower-case, upper-case and calligraphic upper-case letters, respectively. The transpose of a given matrix \textbf{A}, is written as $\mathbf{A}^T$.
An entry of a vector \textbf{a}, a matrix \textbf{A}, or a tensor $\boldsymbol{\mathcal{A}}$ is denoted by $a_i$, $a_{i,j}$,  $a_{i,j,k}$, etc. The symbol $\otimes$ denotes the Kronecker product. The column-wise Khatri-Rao product of two matrices, $\mathbf{A} \in \mathbb{R} ^{I\times R}$ and $\mathbf{B} \in \mathbb{R} ^{J\times R}$, is denoted by $\mathbf{A}\odot\textbf{B}=\begin{bmatrix}\mathbf{a}_1\otimes \mathbf{b}_1, \mathbf{a}_2\otimes \textbf{b}_2,\ldots, \textbf{a}_R\otimes \textbf{b}_R \end{bmatrix}$ with $\textbf{a}_i,\textbf{b}_i$ being the $i$th columns of $\mathbf{A},\mathbf{B}$, respectively. The outer product is denoted by $\circ$. For an $N$th-order tensor, $\boldsymbol{\mathcal{A}} \in \mathbb{R} ^{I_1 \times I_2 \times \cdots \times I_N}$, $\textbf{A}_{\times n}$ denotes its mode−$n$ unfolded (matricized) version resulting from mapping the tensor element with indices $(i_1,i_2,\ldots,i_N)$ to a matrix element $(i_n,j)$ with $j=1 + \sum_{k=1,k\neq n}^N [ ( i_k -1 ) \prod_{m=1,m\neq n}^{k-1}I_m ]$ for $N>2$.

\section{Tensorial fMRI Analysis}
\label{sec:theory}

This section briefly outlines the most commonly used tensorial fMRI methods. Their advantages and disadvantages are discussed, with an emphasis on the reasons that justify the present investigation. Three-way tensors will be considered in this section, for simplicity. 

\subsection{Canonical Polyadic Decomposition (CPD)}
The Canonical Polyadic Decomposition (CPD) (or PARAFAC or CANDECOMP)~\cite{1970_carroll_analysis,1970_harshman_foundations} approximates a tensor, $\boldsymbol{\mathcal{T}} \in \mathbb{R} ^{I_1 \times I_2\times I_3}$, by a sum of $R$ rank-1 tensors (vector outer products),
\begin{equation}
\label{cpd1}
\boldsymbol{\mathcal{T}} = \sum_{r=1}^{R} \mathbf{a}_r \circ \mathbf{b}_r \circ \mathbf{c}_r + \boldsymbol{\mathcal{E}},
\end{equation}
where $\boldsymbol{\mathcal{E}}$ stands for the modeling error tensor. 
Equivalently \cite{2005_smilde_multi-way}, 
\begin{equation}
\label{cpd2}
\textbf{T}_{\times 1} = \mathbf{A} (\mathbf{C}\odot\mathbf{B})^T + \mathbf{E}_{\times 1},
\end{equation}
where $\mathbf{A}:=\begin{bmatrix}\mathbf{a}_1,\mathbf{a}_2,\ldots,\mathbf{a}_R\end{bmatrix}$ and $\mathbf{B},\mathbf{C}$ are similarly defined.
Usually, the sum-of-squares of the residual $\boldsymbol{\mathcal{E}}$ is minimized to determine the latent factors in the $R$ terms.
The main advantage of the CPD, besides its simplicity, is the fact that it is unique up to permutation and scaling under mild conditions. The most widely known and used sufficient uniqueness condition is based on the k-ranks~\cite{1977_kruskal_three-way} of the factor matrices: 
\begin{equation}
\label{eqkrus}
k_\textbf{A}+k_\textbf{B}+k_\textbf{C} \geq 2R +2.
\end{equation}
The uniqueness of CPD is crucial to its application in fMRI analysis. In fact, it was demonstrated (see~\cite{2007_stegeman_comparing,2013_helwig_critique}) that CPD with fMRI data is robust to overlaps (spatial and/or temporal) as well as noise. On the other hand, the result of the CPD method is largely dependent on the correct estimation of the tensor rank $R$~\cite{2003_bro_new,2002_castellanos_triangle,2016_liu_detection,2007_li_estimating}, which is a well known difficult task in tensor modeling. 
 
\subsection{Tensor Probabilistic Independent Component Analysis (TPICA)}

Independent Component Analysis (ICA) is a powerful tool for separating a multivariate signal into additive components, based on the assumption that they are independent. The well known Principal Component Analysis (PCA) is given a stochastic interpretation and the assumption of orthogonal components is strengthened to statistically independent components. ICA has demonstrated promising results in the characterization of fMRI data~\cite{2006_calhoun_unmixing}. TPICA, as proposed in ~\cite{2005_beckmann_tensorial}, is essentially a hybrid of the Probabilistic ICA (PICA)~\cite{2004_beckmann_probabilistic} method and the CPD method. Given a tensor of fMRI data, $\boldsymbol{\mathcal{T}}$, of the form voxels × time × subjects, TPICA factorizes it as:
\begin{equation}
\label{tpica}
\mathbf{T}_{\times 1} \approx \mathbf{A}\mathbf{M}^T,
\end{equation}
where $\mathbf{A}$ is a matrix that contains the weights of the $R$ components of the $I$ voxels (spatial maps), such that each row is assumed to be a sample of $R$ independent, non-Gaussian (assumption following from~\cite{2004_beckmann_probabilistic}) random variables and $\mathbf{M}$ is a Khatri-Rao structured mixing matrix,
\begin{equation}
\mathbf{M} =\mathbf{C}\odot\mathbf{B}.
\label{mixmax}
\end{equation}
TPICA computes the decomposition of the tensor $\boldsymbol{\mathcal{T}}$ in two steps, as dictated by eqs.~(\ref{tpica}) and~(\ref{mixmax}). First, an ICA step is performed, which estimates $\mathbf{M}$ and $\mathbf{A}$, and, at a second step, a Khatri-Rao factorization of $\mathbf{M}$ is computed (using Singular Value Decomposition (SVD)) to compute the factors $\mathbf{B}$ and $\mathbf{C}$.

TPICA is more robust than CPD to rank estimation inaccuracies but it exhibits inferior performance in the presence of overlap in the sources and/or strong noise~\cite{2007_stegeman_comparing,2013_helwig_critique}.

\section{Block Term Decomposition (BTD) for Spatially Folded fMRI Data}
\label{sec:BTD}

Phan et al.~\cite{2013_phan_candecomp/parafac,2013_tichavsky_cramer-rao-induced} proved that, unfolding noisy data to low dimensional tensors (including matrices and vectors as special cases) results in loss of accuracy in the respective decomposition, unless the tensor has a mode with orthogonal components. In fMRI orthogonality of the components is not relevant and hence a loss of accuracy is inevitably, the extent of this loss in accuracy depends on the collinearity of the components of the resulted modes. 

Furthermore, as pointed out in~\cite{1995_norgaard_classification}, the ability of multiway fits to make predictions more robustly, compared to their two-way counterparts, seems to grow with the noise level of the data. This is an important observation with implications for neuroimaging data in general, and fMRI in particular. The use of tensors of higher order could, therefore, improve the result of the decomposition, both in terms of accuracy and robustness in cases of strong noise. 

\begin{figure} [b]
\centering
\captionsetup{justification=centering}
\includegraphics[width=0.99\textwidth]{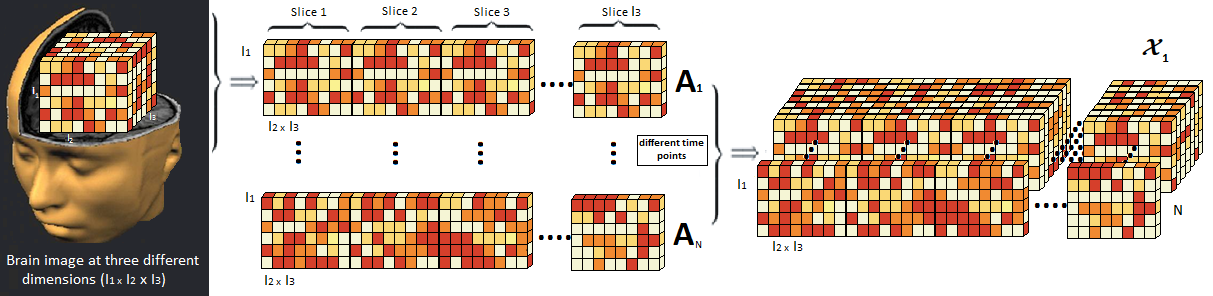}
\caption{Brain images are unfolded in matrices, which are then stacked in tensors. Single subject case} \label{fig:4d}
\end{figure}

A way to benefit from the findings mentioned before is to adopt an alternative type of data unfolding, instead of reshaping the whole brain volume (at a time instance $n$) into a vector $\mathbf{t}_n$ (Fig.~\ref{fig:concat}). Recall that our data, by their very nature, comprise a three-dimensional tensor (even without considering the time dimension). For the proposed unfolding, we adopt the mode-1 (frontal) matricization of the respective data tensor. The different frontal slices are concatenated to form a big matrix $\mathbf{A}_n$ (Fig.~\ref{fig:4d}) (mode-1 matricization of the original spatial tensor). By stacking these matrices together a 3rd-order tensor, $\boldsymbol{\mathcal{X}}_n$, is formed (Fig.~\ref{fig:4d}). Finally, for different subjects, a 4-th order tensor is created by stacking all three-dimensional $\boldsymbol{\mathcal{X}}_n$ tensors. 

The CPD for the four-dimensional data can then be expressed as~\cite{2014_favier_overview}

\begin{figure}
\centering
\captionsetup{justification=centering}
\includegraphics[width=0.99\textwidth]{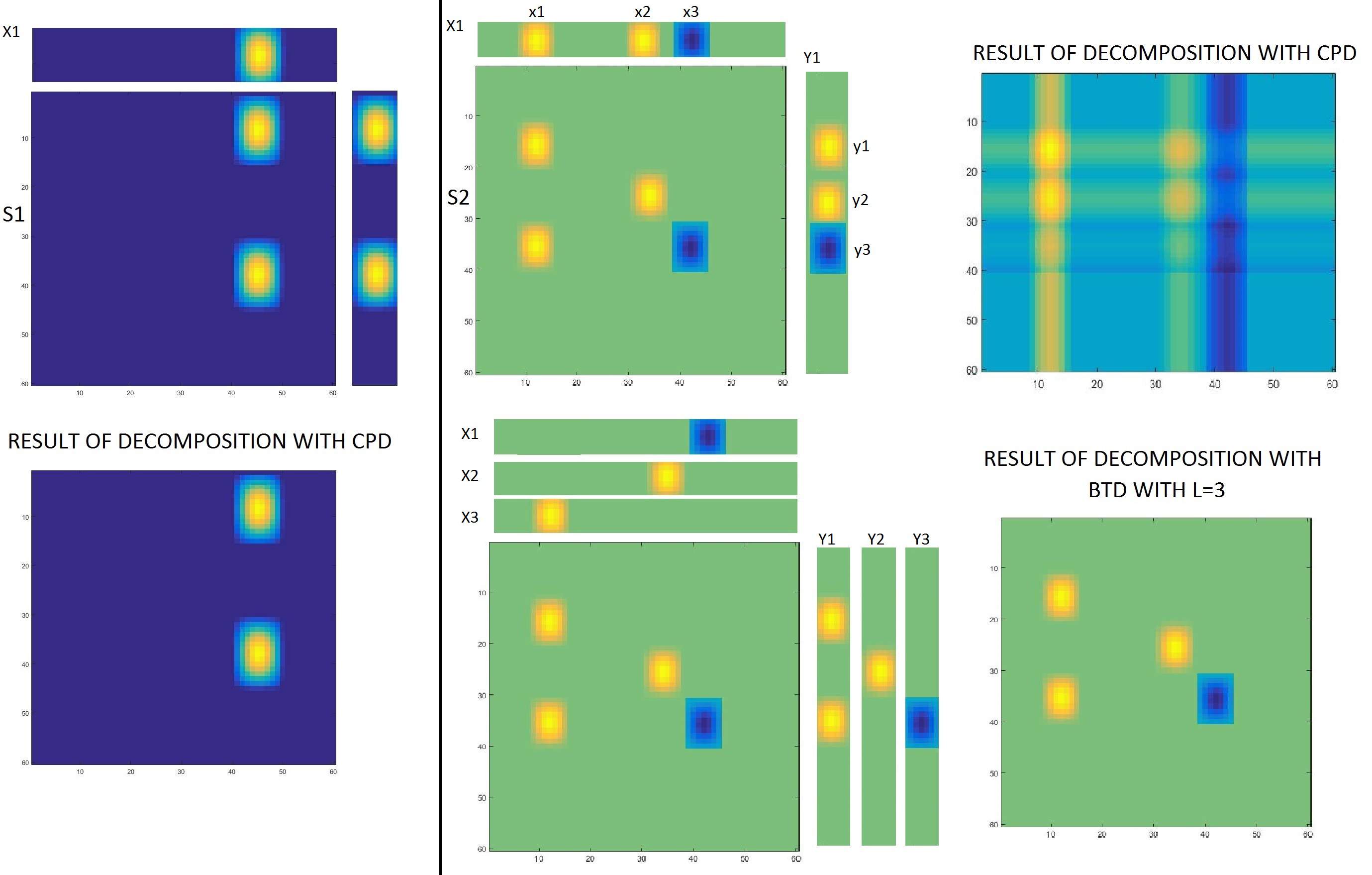}
\caption{\label{fig:examp} Decomposition of sources $S_1$ (left) and $S_2$ (right).}
\end{figure}

\begin{equation}
\label{cpd4d1}
\boldsymbol{\mathcal{T}}=\sum_{r=1}^{R} \mathbf{x}_r \circ \mathbf{y}_r \circ \mathbf{b}_r \circ \mathbf{c}_r + \boldsymbol{\mathcal{E}},
\end{equation}
or equivalently
\begin{equation}
\label{cpd4d2}
\mathbf{T}_{\times 1} = \mathbf{X} (\mathbf{C}\odot \mathbf{B}\odot \mathbf{Y})^T + \mathbf{E}_{\times 1}.
\end{equation}
However, the use of CPD with this specific type of unfolding can be problematic in cases where the spatial maps are not of rank one. An explanation and an example of the problematic nature of the use of CPD on the folded data can be seen in Fig.~\ref{fig:examp}. 
The rank-1 decomposition (CPD) succeeds only in cases where the assumption of rank-1 is correct (S1 source). In cases of spatial activation of higher order (S2 source), the wrong result of applying CPD directly on the folded data can be seen, as several phantoms are produced, from the multiplication of $x_1$ with $y_2$, $x_2$ with $y_1$ and $y_3$, etc. The constraints of CPD can be proved to be restrictive in such cases. 

On the other hand, it seems less restrictive to decompose the tensor in terms of low rank factors, which, however, are not necessarily of rank one; this enhances the potential for modeling more general phenomena~\cite{2012_de_lathauwer_block}.
As an alternative to CPD, the use of Block Term Decomposition (BTD)~\cite{2008_de_lathauwer_decompositions,2008_de_lathauwer_decompositions-1,2008_de_lathauwer_decompositions-2,2011_de_lathauwer_blind} in the folded higher-dimensional tensor is investigated in this work. The rank-($L_r,L_r,1$) BTD of a third-order tensor, $\boldsymbol{\mathcal{T}} \in \mathbb{R} ^{I_1 \times I_2 \times I_3}$, into a sum of rank-($L_r,L_r,1$) terms is given by
\begin{equation}
\label{btd}
\boldsymbol{\mathcal{T}}=\sum_{r=1}^{R} \textbf{A}_r \circ \textbf{b}_r = \sum_{r=1}^{R} (\textbf{X}_r \textbf{Y}_r^T) \circ \textbf{b}_r,
\end{equation}
where the matrix $\textbf{A}_r=\textbf{X}_r \textbf{Y}_r^T \in \mathbb{R} ^{I_1 \times I_2}$ has rank $L_r$. BTD has been successfully applied in modeling epileptic seizures in electro-encephalograms (EEG)~\cite{2014_hunyadi_block} and proved capable of modeling nonstationary seizures, which evolve either in frequency or in space, better than other techniques. BTD has not been previously applied in fMRI analysis, to the best of the authors' knowledge. 

In the common CPD-based method, once the different factors in~(\ref{cpd1}) have been estimated, the spatial mode $\textbf{a}_r$ is reshaped, folded back from a vector to a tensor; in this way, activation patterns in the different frontal slices of the brain are formed. Consider, for example, fMRI images of $60\times 60\times 40=144000$ voxels with duration equal to 100~seconds, in a multi-subject analysis involving 10~different subjects. The ``traditional'' approach would suggest to decompose the tensor $\boldsymbol{\mathcal{T}} \in \mathbb{R} ^{144000 \times 100 \times 10}$ into $R$ vectors $\textbf{a}_r\in \mathbb{R} ^{144000}$, $\textbf{b}_r\in \mathbb{R} ^{100}$ and $\textbf{c}_r\in \mathbb{R} ^{10}$ and then fold each vector $\textbf{a}_r$ back to a tensor $\boldsymbol{\mathcal{A}}_r\in \mathbb{R} ^{60\times 60 \times 40}$. In this work, it is suggested to decompose directly the higher-dimensional data, once the frontal slices have been stacked one after another (i.e., unfolding only one of the three spatial dimensions) instead of first unfolding and then folding back to the original dimensions. The corresponding rank-($L_r,L_r,1,1$) BTD~(\ref{btd}) can be written as
\begin{equation}
\label{btd4d}
\boldsymbol{\mathcal{T}}=\sum_{r=1}^{R} \mathbf{A}_r \circ \mathbf{b}_r \circ \mathbf{c}_r = \sum_{r=1}^{R} (\mathbf{X}_r {\mathbf{Y}_r^T}) \circ \mathbf{b}_r \circ \mathbf{c}_r,
\end{equation}
where the matrix $\mathbf{A}_r=\mathbf{X}_r {\mathbf{Y}_r^T} \in \mathbb{R} ^{60 \times 2400}$ has rank $L_r$. 

\subsection{Uniqueness}

It was proved in~\cite{2008_de_lathauwer_decompositions-2}
that the BTD~(\ref{btd}) is essentially unique up to scaling, permutation and the simultaneous post-multiplication of $\textbf{X}_r$ by a nonsingular matrix $\textbf{F}$ with the pre-multiplication of $\textbf{X}_r$ by $\textbf{F}^{-1}$ provided that the matrices $\begin{bmatrix}\begin{array}{cccc} \textbf{X}_1 & \textbf{X}_2 & \cdots & \textbf{X}_R \end{array}\end{bmatrix}$ and $\begin{bmatrix}\begin{array}{cccc} \textbf{Y}_1 & \textbf{Y}_2 & \cdots & \textbf{Y}_R \end{array}\end{bmatrix}$ are full column rank and the matrix $\mathbf{B}=\begin{bmatrix}\begin{array}{cccc}\textbf{b}_1 & \textbf{b}_2 & \cdots & \textbf{b}_R\end{array}\end{bmatrix}$ does not contain collinear columns. A simple argument showing that this uniqueness condition can be extended to the rank-$(L_r,L_r,1,1)$ case follows.

\begin{proof}[Proof of uniqueness of rank-$(L_r,L_r,1,1)$ BTD]
As suggested in~\cite{2000_sidiropoulos_uniqueness}, the uniqueness of a higher-order tensor decomposition can be shown through a reduction to a third-order tensor, which is ``the first instance of multilinearity for which uniqueness holds and from which uniqueness propagates by virtue of Khatri-Rao structure'' \cite{2000_sidiropoulos_uniqueness}. 

Assume that the matrices $\left[\begin{array}{cccc} \textbf{X}_1 & \textbf{X}_2 & \cdots & \textbf{X}_R\end{array}\right]$ and $\left[\begin{array}{cccc} \textbf{Y}_1 & \textbf{Y}_2 & \cdots & \textbf{Y}_R\end{array}\right]$ are of full column rank and the matrices $\textbf{B}$, $\mathbf{C}=\begin{bmatrix}\begin{array}{cccc}\textbf{c}_1 & \textbf{c}_2 & \cdots & \textbf{c}_R\end{array}\end{bmatrix}$ do not contain collinear or null columns (a realistic assumption for matrices that represent time and subjects). In view of the above, uniqueness for the unfolded version of~(\ref{btd4d}), which has the fourth mode nested into the third one, must be proved:
\begin{equation}
\label{equn1}
\boldsymbol{\mathcal{T}}=\sum_{r=1}^{R} \textbf{A}_r \circ  \textbf{g}_r,
\end{equation}
where
\begin{equation}
\label{equn2}
\textbf{G}=\left[\begin{array}{cccc} \textbf{g}_1 & \textbf{g}_2 & \cdots & \textbf{g}_R \end{array}\right]=\textbf{B} \odot \textbf{C}.
\end{equation}
Eq.~(\ref{equn1}) is the BTD of a three-way tensor as in~(\ref{btd}). The matrix $\textbf{G}$, equal to the Khatri-Rao product of two matrices, which do not have null or collinear columns, does not have null or collinear columns either (see, e.g.,~\cite[Proposition 1]{2011_brie_uniqueness}).

Following Theorem~4.1 of~\cite{2008_de_lathauwer_decompositions-2}, since the matrices $\left[\begin{array}{cccc} \textbf{X}_1 & \textbf{X}_2 & \cdots & \textbf{X}_R\end{array}\right]$ and $\left[\begin{array}{cccc} \textbf{Y}_1 & \textbf{Y}_2 & \cdots & \textbf{Y}_R\end{array}\right]$ have full column rank and the matrix $\textbf{G}$ has no collinear columns, the decomposition is unique.
\end{proof}

\section{Simulation Results}

Three different simulation studies are presented, with scenarios and data reproduced from~\cite{2013_helwig_critique,2007_stegeman_comparing,2011_erhardt_comparison}. The first experiment~\cite{2013_helwig_critique} simulates a perception study in a scenario with a single-slice brain and consists of three different sources (with temporal and spatial overlap) and Gaussian noise. The experiment in~\cite{2007_stegeman_comparing} (adapted from the experiment in~\cite{2005_beckmann_tensorial}) consists of three different sources either without overlap or with significant overlap, in three different slices (so that crosstalk of the spatial maps among different slices can be checked) with the presence of resting state fMRI noise and Gaussian noise. Finally, the experiment from~\cite{2011_erhardt_comparison} simulates (as in SimTB~\cite{2012_erhardt_simtb}) eight sources with different types of noise  but not significant overlap in the spatial maps.

Rank determination in the experiments is performed with the aid of the Core Consistency Diagnostic (CorConDia) method~\cite{2003_bro_new} and the triangle method~\cite{2002_castellanos_triangle}, implemented by Tensorlab~3.0~\cite{2016_vervliet_tensorlab}. In all the experiments, the estimated rank of the decomposition depends on the noise level. It increases as the signal to noise ratio (SNR) levels decrease significantly, because some peaks of the noise get higher amplitude than the useful signal and are recognized as a source~\cite{2007_stegeman_comparing}. For the computation of the BTD, the Non Linear Least Squares method of the Structured Data Fusion (SDF) toolbox~\cite{2015_sorber_structured} is employed, as implemented in Tensorlab~3.0~\cite{2016_vervliet_tensorlab}.

\subsection{Simulation of a perception study}

In this study, the data were simulated under the following assumptions: (a) three subjects have participated in this simplified version of a realistic perception study~\cite{2013_helwig_critique} lasting 180~seconds. The subjects watched a gray dot appearing on a display screen for 30-sec periods, followed by 30-sec periods of a blank screen; (b) the second (of the three) displayed dots moved on the screen during its 30-sec period, whereas the first and third displayed dots were stationary on the screen; (c) the subjects were instructed to verbally indicate when the dot has disappeared from the screen. 

\begin{figure}
\centering
\captionsetup{justification=centering}
\includegraphics[width=0.69\textwidth]{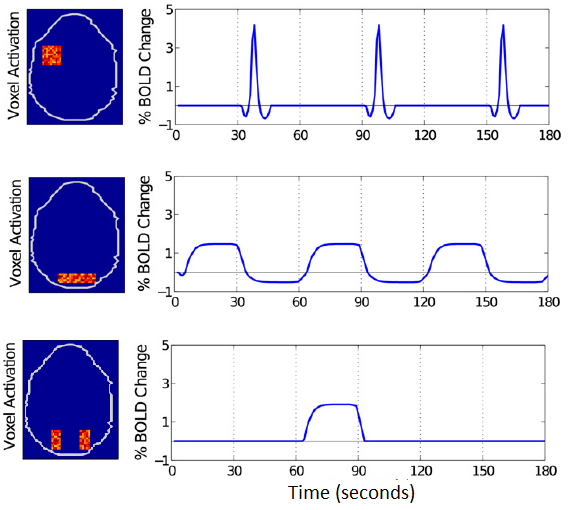}
\caption{\label{fig:sources} The three sources used in the first experiment.}
\end{figure}

The simulated data are a $60\times 50$ axial slice of voxel activity (3000~voxels in total) from somewhere near the level of Broca’s area. The second and third components of Fig.~\ref{fig:sources} represent the visual and motion perception components respectively, which have $50\%$ of shared active voxels. The data from each subject contained all the three sources with different activation levels equal to

\begin{equation}
\label{activation}
\textbf{C}=\left[\begin{array}{cccc} 3 & 4 & 5 \\ 2 & 3 & 4 \\2 & 2 & 3 \end{array}\right],
\end{equation}

\noindent times the mean noise standard deviation. The activity within each active voxel was randomly sampled from a Uniform~[0.8,1.2] distribution for each replication of each simulation condition. Finally, note that approximately $3$-$6\%$ of the voxels are active at any point during the experiment. Gaussian noise was added in all cases. Following~\cite{2007_stegeman_comparing}, the SNR is defined as the Frobenius norm of the signal divided by the Frobenius norm of the noise:
\begin{equation}
\label{snr}
\mathrm{SNR}=\frac{\|\textbf{A} (\textbf{C}\odot\textbf{B})^T\|}{\|\textbf{E}\|}.
\end{equation}

\begin{figure}
\centering
\captionsetup{justification=centering}
\includegraphics[width=0.99\textwidth]{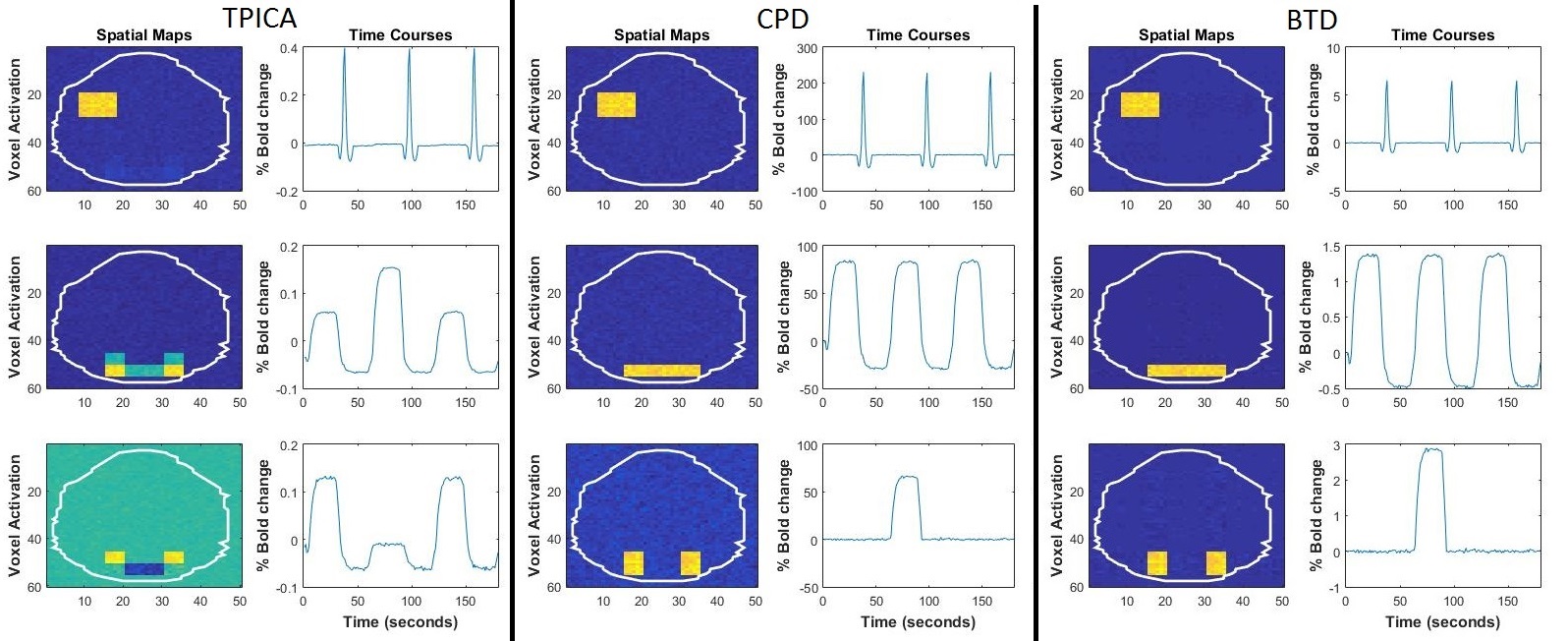}
\caption{\label{fig:1.2}  Decomposition of the data (SNR=0.12).}
\end{figure}

\begin{figure}
\centering
\captionsetup{justification=centering}
\includegraphics[width=0.99\textwidth]{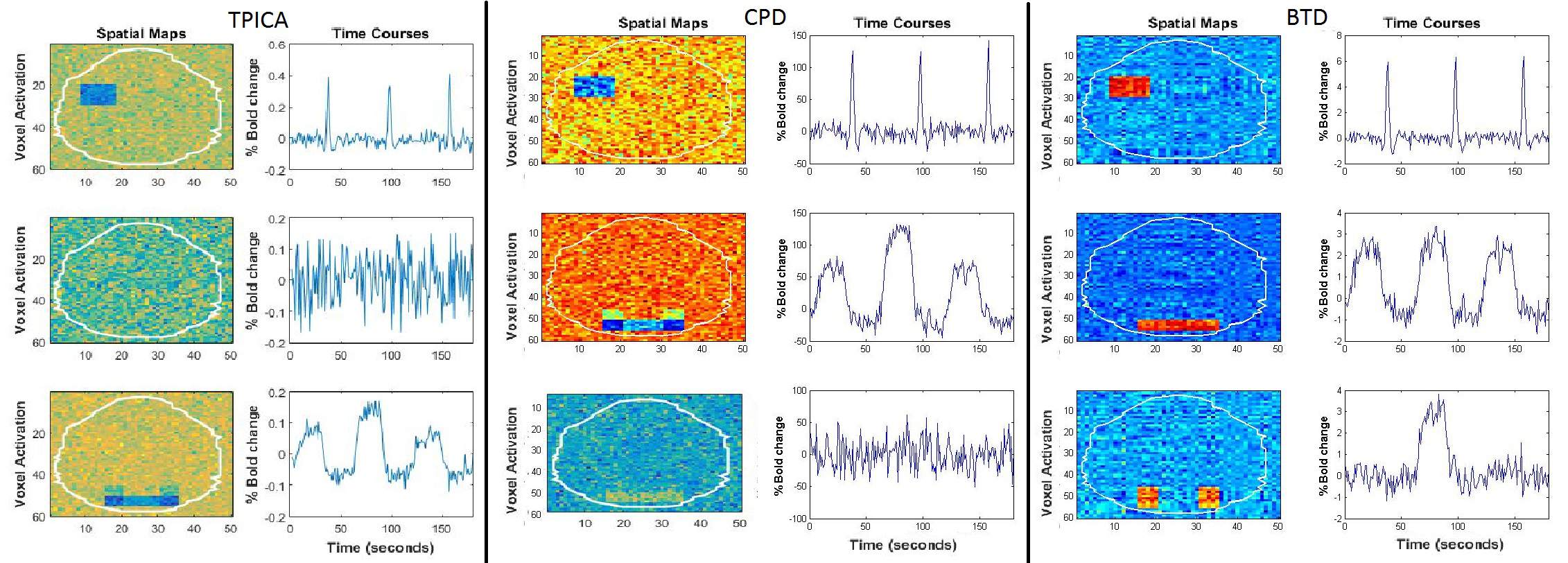}
\caption{\label{fig:0.8}  Decomposition of the data (SNR=0.08).}
\end{figure}

\begin{figure}
\centering
\captionsetup{justification=centering}
\includegraphics[width=0.99\textwidth]{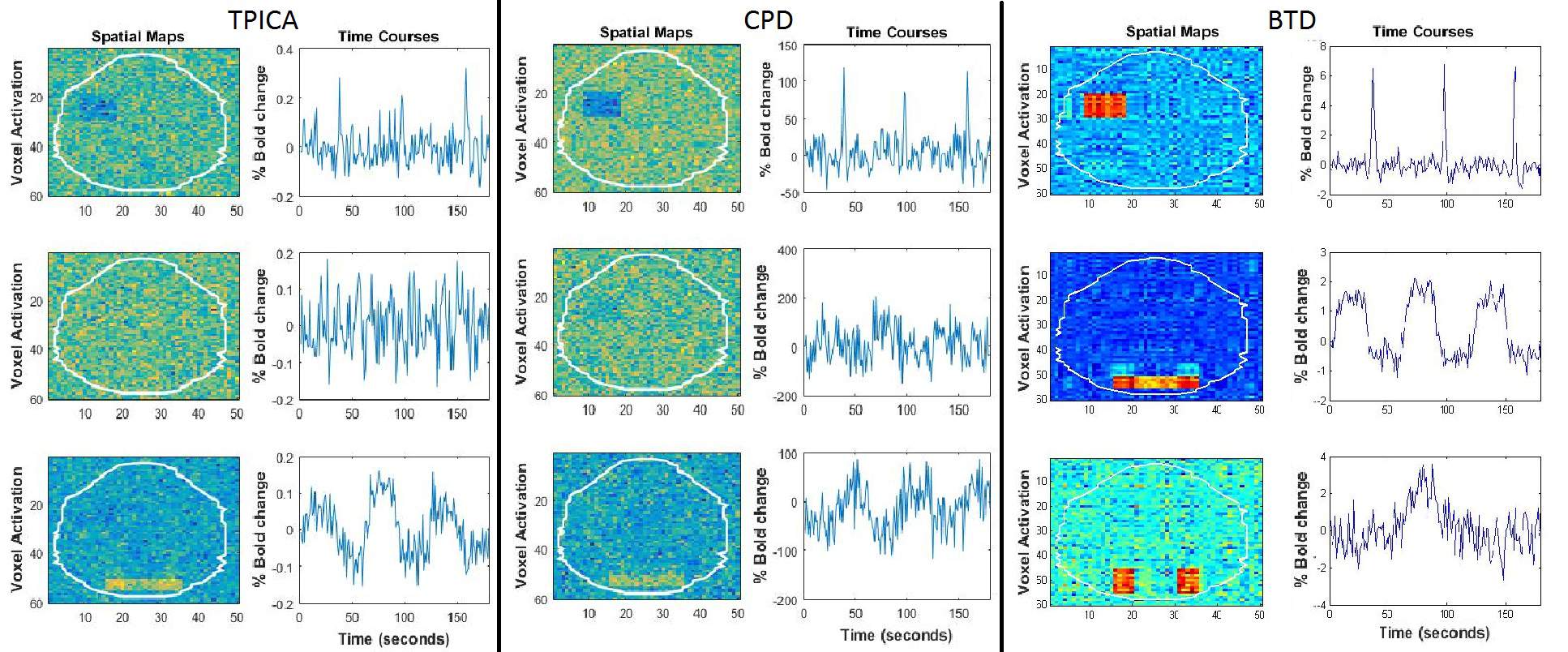}
\caption{\label{fig:0.5}  Decomposition of the data (SNR=0.05).}
\end{figure}
With all the methods, three components (rank of decomposition) were extracted. The rank $L_r$ of the BTD is set equal to three (for all $r$). The method is not sensitive to overestimation of $L_r$; as a matter of fact, one gets same (when the true rank of the source is equal to or smaller than $L_r$) or better (when the true rank of the source is much larger than $L_r$) results with higher values of $L_r$, at the cost of increased complexity. A compromise between accuracy and complexity shall be made. In Figs.~4-6, the performance gain of BTD as the noise power increases can be observed. TPICA starts failing at values of SNR lower than 0.12 while CPD, following the findings of \cite{2013_helwig_critique}, can maintain the almost perfect separation of the sources at values of $\mathrm{SNR}\approx 0.1$. The fact that the proposed BTD identifies the sources correctly at the overlapping areas even with values of SNR lower by $100\%$ must be emphasized.

\subsection{Multi-slice simulation}

The signal here consists of artificial voxel activation maps, time patterns and activation strengths for three subjects. Each activation map consists of three different slices. Random Gaussian noise is added, and for each voxel the noise mean and variance are estimated from real single-subject resting state fMRI measurements (for details, see~\cite{2005_beckmann_tensorial}). The voxel-wise noise mean and variance are the same for each of the three subjects.

\begin{figure}
  \centering
  \captionsetup{justification=centering}
  \includegraphics[width=.85\linewidth]{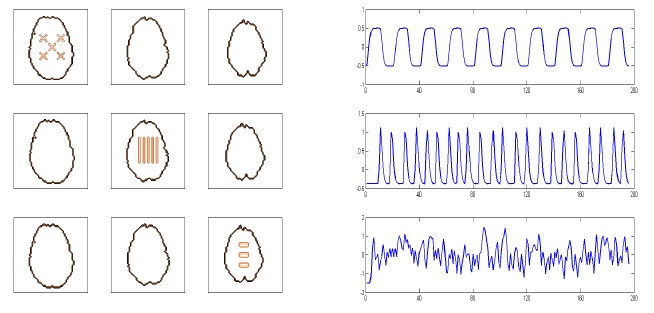}
  \caption [width=.79\linewidth]{ Sources of dataset A \cite{2005_beckmann_tensorial}.}
  \label{fig:stega}
\end{figure}

Beckmann and Smith~\cite{2005_beckmann_tensorial} consider five different artificial fMRI datasets, named (A)-(E), which differ only in their signal part, while A. Stegeman~\cite{2007_stegeman_comparing} added three more  fMRI datasets (F)-(H) with high percentage of overlap between the sources (datasets used in~\cite{2005_beckmann_tensorial} do not have overlap). In this paper, datasets (A) (Fig.~\ref{fig:stega}) and dataset (G)  (Fig.~\ref{fig:stegg}) will be used. Each spatial map has a different associated time course: spatial map 1 is modulated by a time course of a simple block design, spatial map 2 by a time course which corresponds to a single-event design, with fixed interval, and spatial map 3 by a single-event design where a random interval between the different events has been introduced. After being convolved with a canonical haemodymanic response function, the time courses and the spatial maps consist of three different spatiotemporal processes, which are present in every subject with different power. 

The performance evaluation was based on the Pearson correlation. We computed the correlation between the time courses obtained from the different methods and the actual ones, and the correlation between the spatial maps acquired and the real ``averaged'' spatial maps for every subject computed with Ordinary Least Squares (OLS) regression (similarly to~\cite{2005_beckmann_tensorial,2007_stegeman_comparing}).

\subsubsection{Dataset A}

The activation levels $\mathbf{C}$ for each subject for the different sources (Fig.~\ref{fig:stega}) are the same as in~(\ref{activation}). The number of active voxels in each of the three spatial maps is 45 out of 962, 90 out of 838 and 54 out of 689, respectively. Hence, only approximately eight percent of all voxels is active, which is a relatively low percentage. In all the methods four components were extracted, due to the fact that this is the best rank for CPD (as exhibited in \cite{2007_stegeman_comparing}, due to the fact that one noise peak is considered as signal because is higher than one of the fMRI signals) and TPICA is robust to different ranks.

Tables~\ref{tab:da} and~\ref{tab:ta} present the mean correlation (over 10 runs) of spatial maps and time courses (respectively) of the different decompositions at three different SNR values. Note the stability in the performance of BTD compared to the other two methods. Furthermore, the different effect of noise in TPICA and CPD can be observed. In TPICA, the correlation between the estimated spatial map and the ``true'' one decreases dramatically as the level of the noise gets higher, while in CPD the decrease is lower but with a significant increase of the cross-talk (correlation between ``wrong'' spatial maps). The correlations compared to those computed in \cite{2007_stegeman_comparing} are slightly higher, since not only the intra-cranial voxels were used but all the voxels inside the minimum rectangle around the brain mask (due to the need for the same number of voxels in every row for BTD). 

\begin{figure} 
  \centering
  \captionsetup{justification=centering}
  \includegraphics[width=1\linewidth]{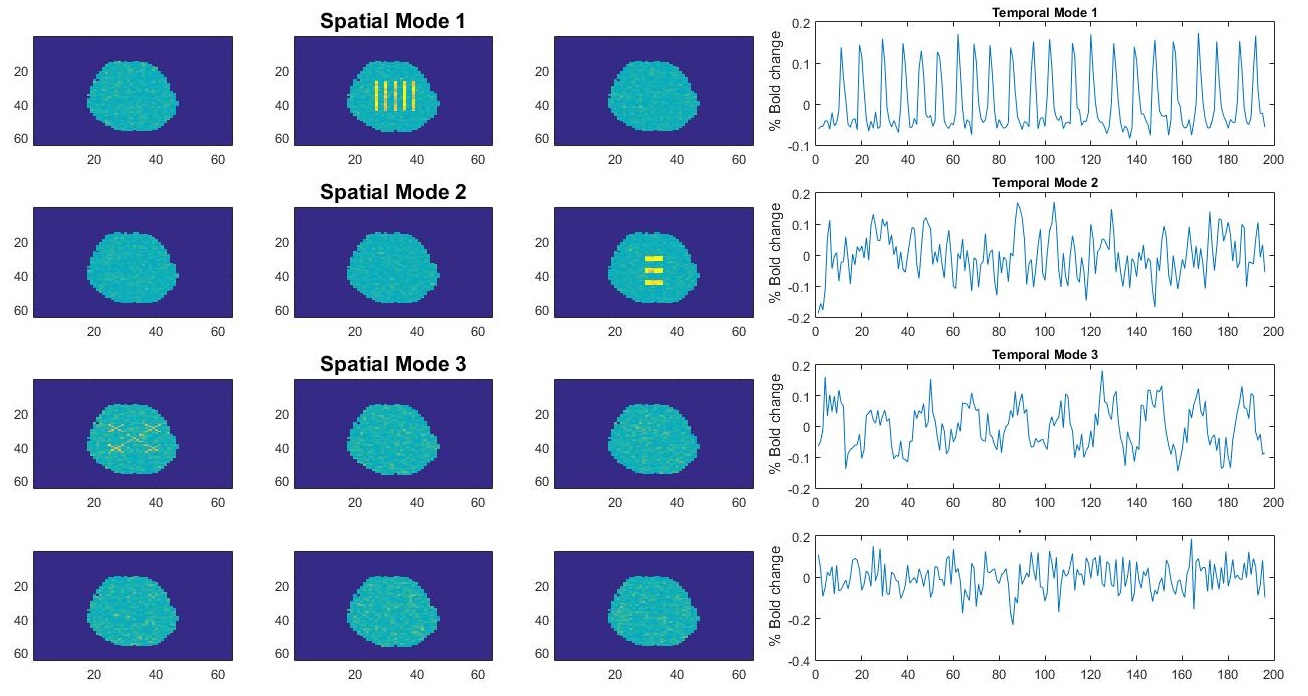}
  \caption [width=.80\linewidth]{TPICA of dataset A with SNR equal to 0.08.}
  \label{fig:tpicA8}
\end{figure}

\begin{figure}
  \centering
  \captionsetup{justification=centering}
  \includegraphics[width=.999\linewidth]{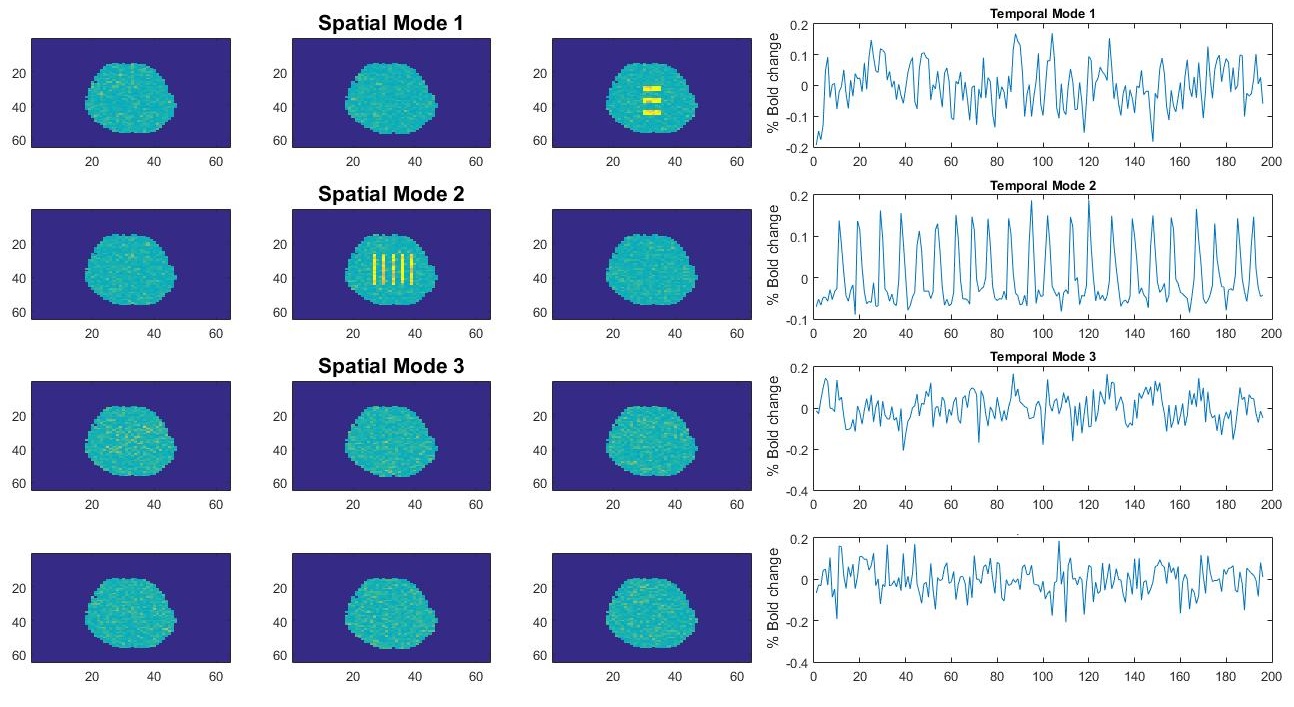}
  \caption [width=.80\linewidth]{TPICA of dataset A with SNR equal to 0.04.}
  \label{fig:tpicA4}
\end{figure}

\begin{figure} 
  \centering
  \captionsetup{justification=centering}
  \includegraphics[width=.999\linewidth]{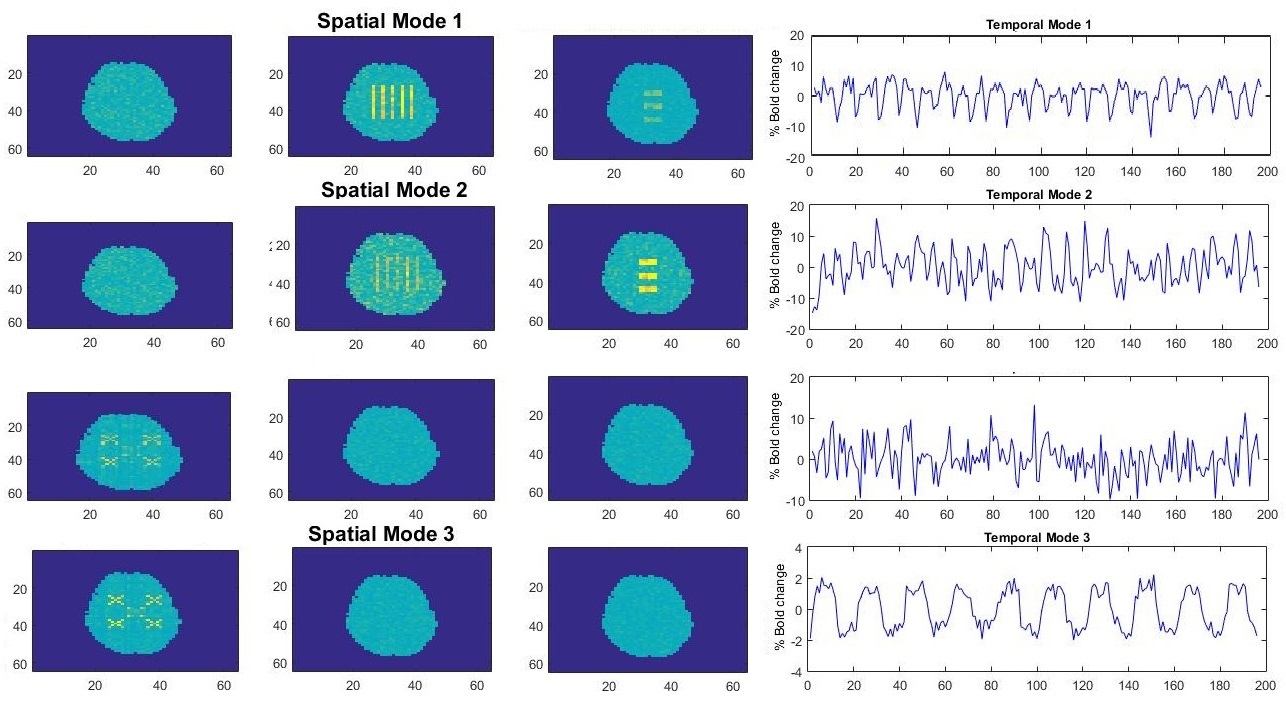}
  \caption [width=.80\linewidth]{CPD of dataset A with SNR equal to 0.08.}
  \label{fig:cpdA8}
\end{figure}

\begin{figure} 
  \centering
    \captionsetup{justification=centering}
  \includegraphics[width=.999\linewidth]{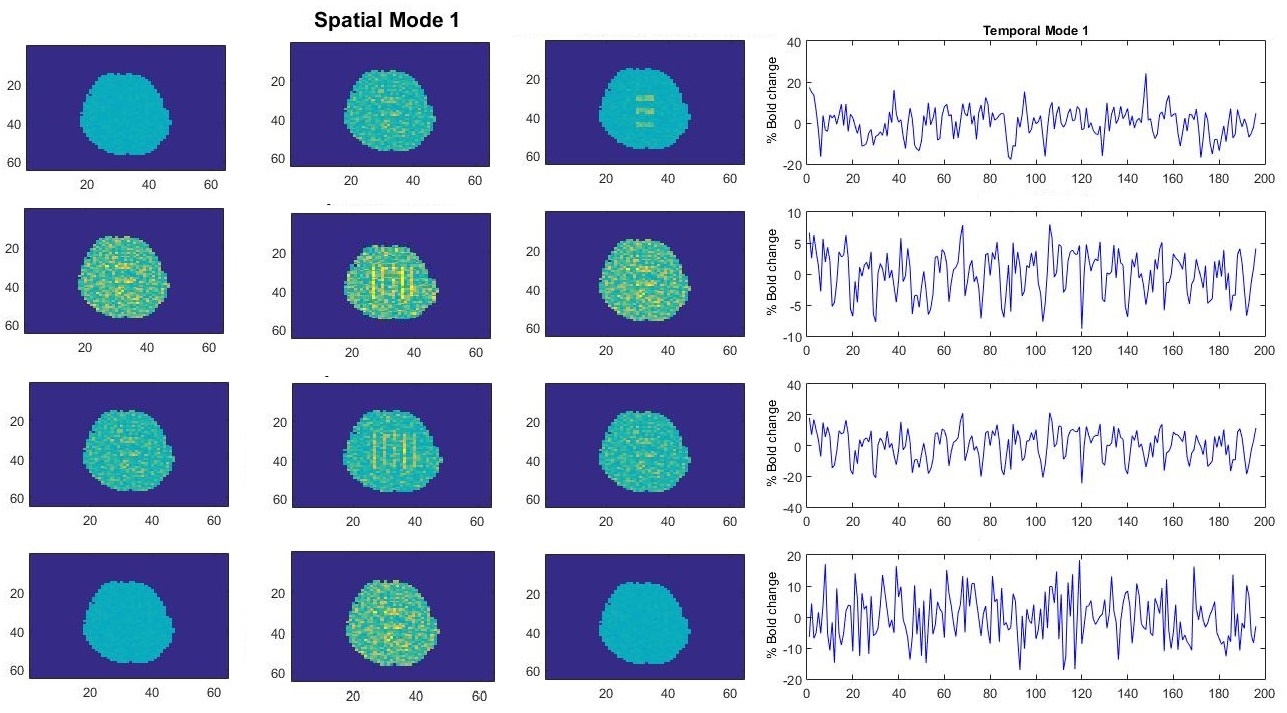}
  \caption [width=.80\linewidth]{CPD of dataset A with SNR equal to 0.04.}
\label{fig:cpdA4}
\end{figure}

\begin{figure} 
  \centering
  \captionsetup{justification=centering}
  \includegraphics[width=.999\linewidth]{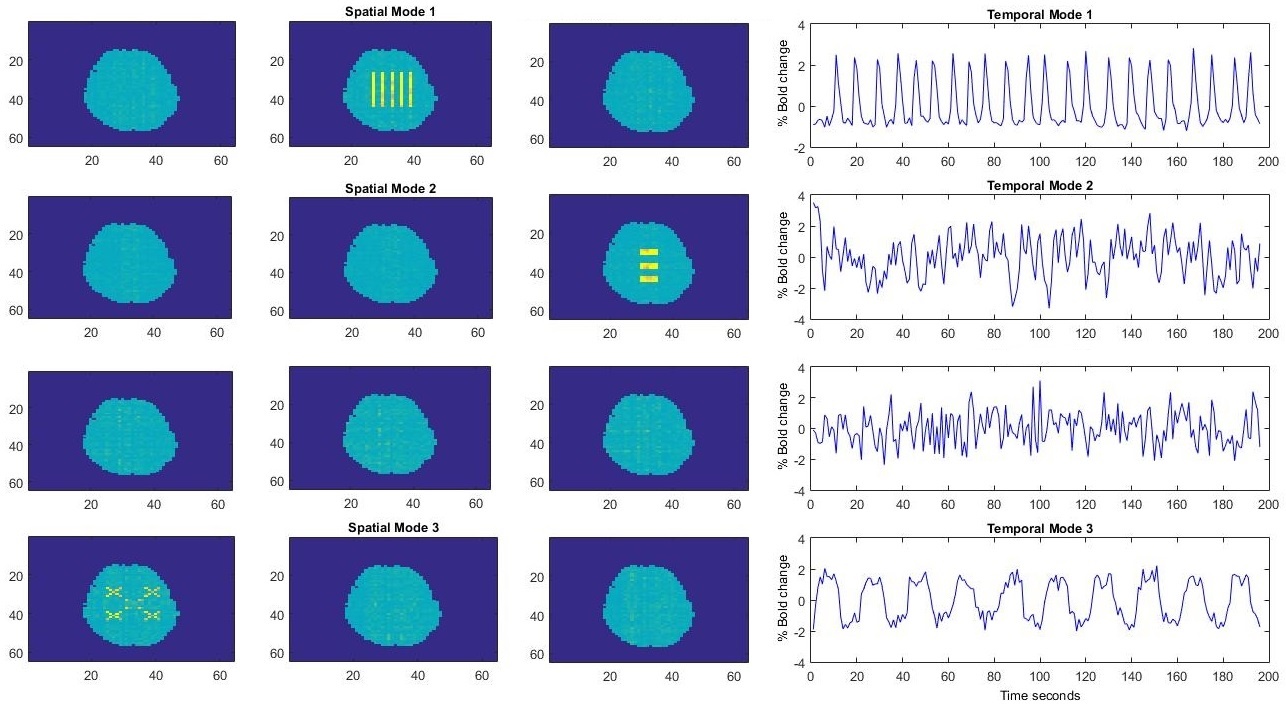}
  \caption [width=.80\linewidth]{BTD of dataset A with SNR equal to 0.08.}
  \label{fig:btdA8}
\end{figure}

\begin{figure} 
  \centering
  \captionsetup{justification=centering}
  \includegraphics[width=.999\linewidth]{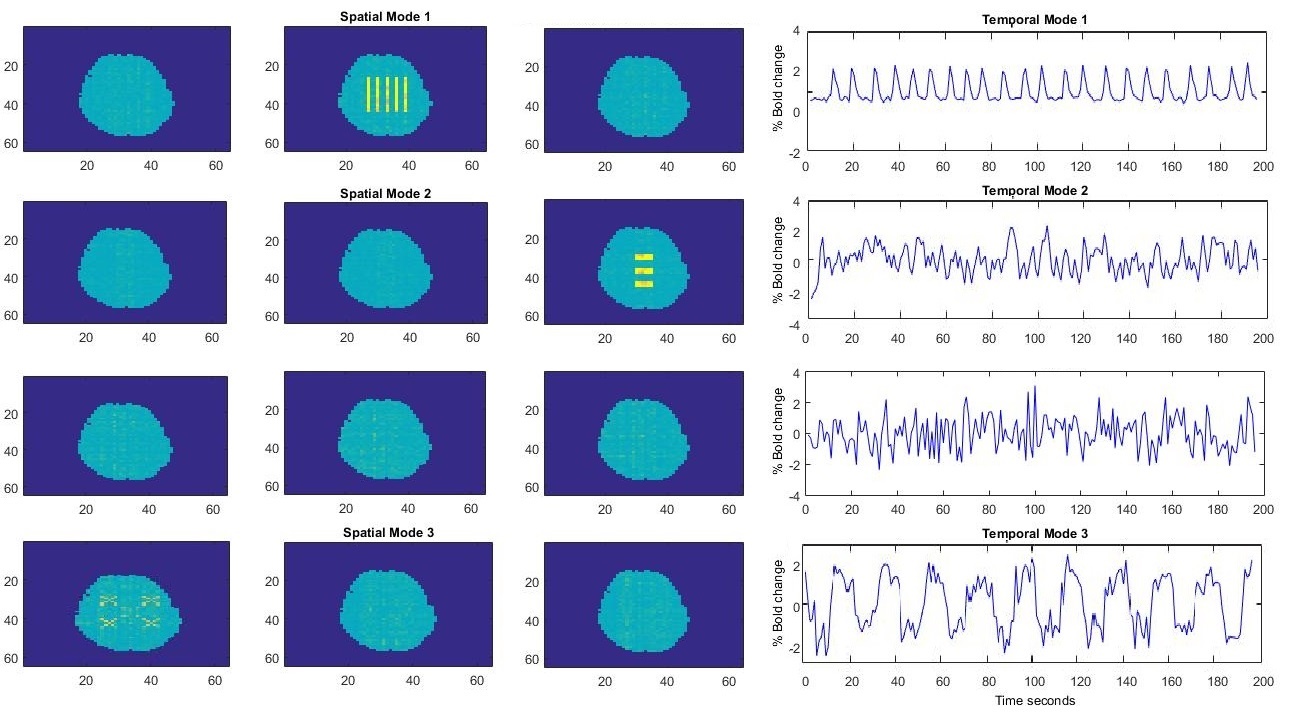}
  \caption [width=.80\linewidth]{BTD of dataset A with SNR equal to 0.04.}
  \label{fig:btdA4}
\end{figure}

Similar conclusions can be reached from Figs.~\ref{fig:tpicA8}--\ref{fig:btdA4}, which present the results of the decompositions at two different SNR values, with the three different methods (the figures depict one of 10 different runs of the algorithms). In every figure, the three different spatial maps and time courses are labeled accordingly, the noise instance does not bear any label. CPD leads to worst results in low noise compared to those of BTD and TPICA with high percentage of cross-talk (as mentioned before) between the different spatial maps; it can be noted that in Fig. ~\ref{fig:cpdA4} only mode 1 can be recognized relatively well. BTD and TPICA have similar performance at SNR=0.08, while at SNR=0.04 the gain in the performance of BTD can be noted (especially in time course 3). It should be mentioned that in this dataset no spatial overlap exists.

  \begin{table} 
  \centering
  \Huge
  \begin{adjustbox}{width=0.8\linewidth}
  \begin{tabular}{||c||c|c|c||c|c|c||c|c|c||}
  \hline
   \multirow{2}{*}{Maps} & \multicolumn{3}{c}{SNR=0.08} & \multicolumn{3}{|c|}{SNR=0.06} & \multicolumn{3}{|c|}{SNR=0.04*} \\
    \hhline{~---------}
    & Map 1 & Map 2 & Map 3 & Map 1 & Map 2 & Map 3 & Map 1 & Map 2 & Map 3  \\
  \hline
  \hline
   \centering TPICA M1 & 0.68 & 0.1 & 0.1 & 0.54 & -0.12 &
  0.15 & 0.34 & 0.19 & 0.18\\ 
  TPICA M2 & 0.1 & 0.99 & 0.1 & 0.12 & 0.89 &
  0.1 & 0.22 & 0.76 & 0.18 \\ 
  TPICA M3 & 0.1 & 0.1 & 0.96& 0.11 & 0.11 &
 0.82 & 0.19 & 0.2 & 0.72 \\
  \hline
     \centering CPD M1 & 0.74 & -0.2 & -0.12 & 0.69 & -0.34& -0.28 & 0.48 & -0.41**& 0.30\\ 
  CPD M2 & 0.18 & 0.97 &  -0.12 & -0.28 & 0.9 & -0.15 & -0.34 & 0.68 & 0.19\\ 
  CPD M3 & -0.21 & -0.12 & 0.99 & -0.29 & 0.15 &
  0.88 & 0.4 & -0.2 & 0.7 \\
  \hline  
  \centering BTD M1 & 0.90 & 0.12 & -0.14 & 0.80 & 0.14 &-0.15 & 0.68 & -0.29 & -0.18\\ 
  BTD M2 & 0.11 & 0.96 & 0.12 & 0.18 & 0.92 & -0.12 & -0.29 & 0.87 & -0.19\\ 
  BTD M3 & -0.12 & 0.1 & 0.99 & -0.2 & 0.12 & 0.94 & -0.23 & -0.19 & 0.84\\
  \hline
  \multicolumn{10}{l}{*At SNR=0.04 four components are extracted}\\
  \multicolumn{10}{l}{**At the fourth component the crosstalk is also relatively high}
\end{tabular}
\end{adjustbox}
  \caption{Correlation of spatial maps of dataset A.}
  \label{tab:da}
\end{table}

  \begin{table} [ht]
  \centering
  \Huge
  \begin{adjustbox}{width=0.8\linewidth}
  \begin{tabular}{||c||c|c|c||c|c|c||c|c|c||}
  \hline
   \multirow{2}{*}{Time Courses} & \multicolumn{3}{c}{SNR=0.08} & \multicolumn{3}{|c|}{SNR=0.06} & \multicolumn{3}{|c|}{SNR=0.04} \\
    \hhline{~---------}
    & TC 1 & TC 2 & TC 3 & TC 1 & TC 2 & TC 3 & TC 1 & TC 2 & TC 3  \\
  \hline
  \hline
\centering TPICA T1 & 0.34 & * & * & 0.28 & * &
  0.1 & 0.2 & * & 0.22\\   
  TPICA T2 & * & 0.92 & * & * & 0.89 &
  *0.12 & * & 0.70 & 0.16 \\ 
  TPICA T3 & 0.1 & 0.1 & 0.96& 0.11 & 0.11 &
 0.82 & 0.19 & 0.2 & 0.72 \\
  \hline
     \centering CPD T1 & 0.58 & * & 0.18 & 0.5 & 0.15 & -0.2 & 0.4 & 0.2 & 0.2\\ 
  CPD T2 & 0.20 & 0.92 & -0.24 & 0.34 & 0.88 & 0.28 & 0.35 & 0.75 & 0.39 \\ 
  CPD T3 & 0.18 & * & 0.95 & -0.2 & 0.12 & 0.85 & 0.22 & -0.23 & 0.68\\
  \hline  
  \centering BTD T1 & 0.68 & * & * & 0.65 & 0.1 & 0.12 & 0.58 & -0.21 & 0.22\\ 
  BTD T2 & 0.11 & 0.92 & 0.18 & -0.24 & 0.89 & 0.22 & 0.23 & 0.85 & 0.28 \\ 
  BTD T3 & -0.18 & * & 0.95 & -0.2 & 0.1 & 0.88 & -0.22 & -0.18 & 0.78\\
  \hline
\end{tabular}
\end{adjustbox}
  \caption{Correlation of time courses of  dataset A.}
  \label{tab:ta}
\end{table}

\clearpage
\subsubsection{Dataset G}

\begin{figure}
  \centering
  \captionsetup{justification=centering}
  \includegraphics[width=.8\linewidth]{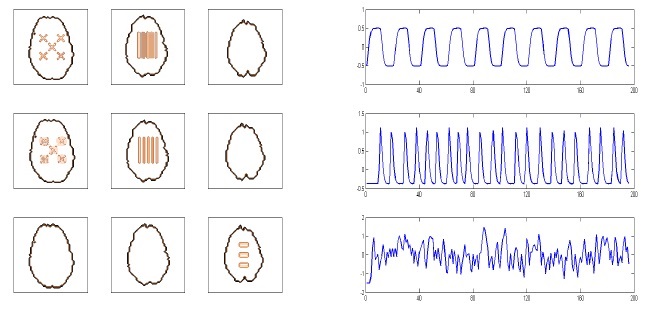}
  \caption [width=.95\linewidth]{ Sources of dataset G.}
  \label{fig:stegg}
  \vspace{-2mm}
\end{figure}

In dataset G, the first spatial activity map is the sum of map~1 in dataset~A and map~2 expanded to the $18\%$ of the active dataset (Fig.~\ref{fig:stegg}). The second activity map is the sum of map~2 in dataset~A and map~1 expanded to the  $18\%$ of the active dataset. Hence, the activity of the first two maps takes place in the first two brain slabs, and they have 51\% and 63\%  of their active voxels in common (high percentage of overlap), respectively. In map 3, the time courses $\mathbf{B}$, the subject levels $\mathbf{C}$ for each subject and the noise instance are the same as in dataset~A.

Tables~\ref{tab:dg}, \ref{tab:tg} exhibit the drawback of CPD and BTD compared to TPICA in cases of high overlap, and also the gain obtained by BTD against any other method, when high levels of noise are present. Even in the first decomposition with relatively high SNR (=0.12), TPICA fails to correctly separate the sources, as the second slice of the spatial maps of the first two sources is almost the same. At SNR=0.08, the maximum correlation of map 2 occurs with spatial map 1, instead of the correct spatial map 2. This is due to the fact that TPICA fails in distinguishing the two different sources, and source 2 is identified as noise (mixture of all sources). On the other hand, CPD and BTD recognize correctly all the sources at high SNR. At lower SNR, CPD can not distinguish the second slice of the two first sources, while BTD still gives an almost perfect result. As noise level increases, the difference between  the results of the decompositions increases in favor of BTD.

\begin{figure}
  \centering
  \captionsetup{justification=centering}
  \includegraphics[width=1\linewidth]{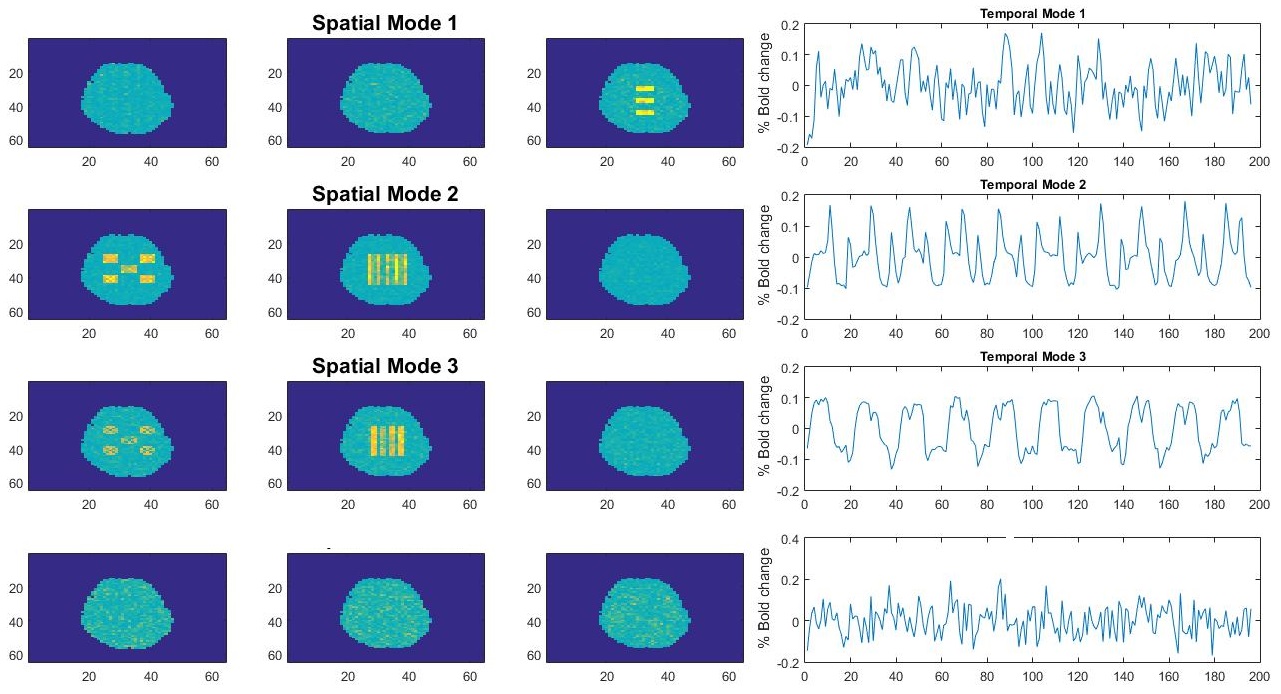}
  \caption [width=.80\linewidth]{TPICA of dataset G with SNR equal to 0.15.}
  \label{fig:tpicG15}
\end{figure}

\begin{figure}
  \centering
  \captionsetup{justification=centering}
  \includegraphics[width=1\linewidth]{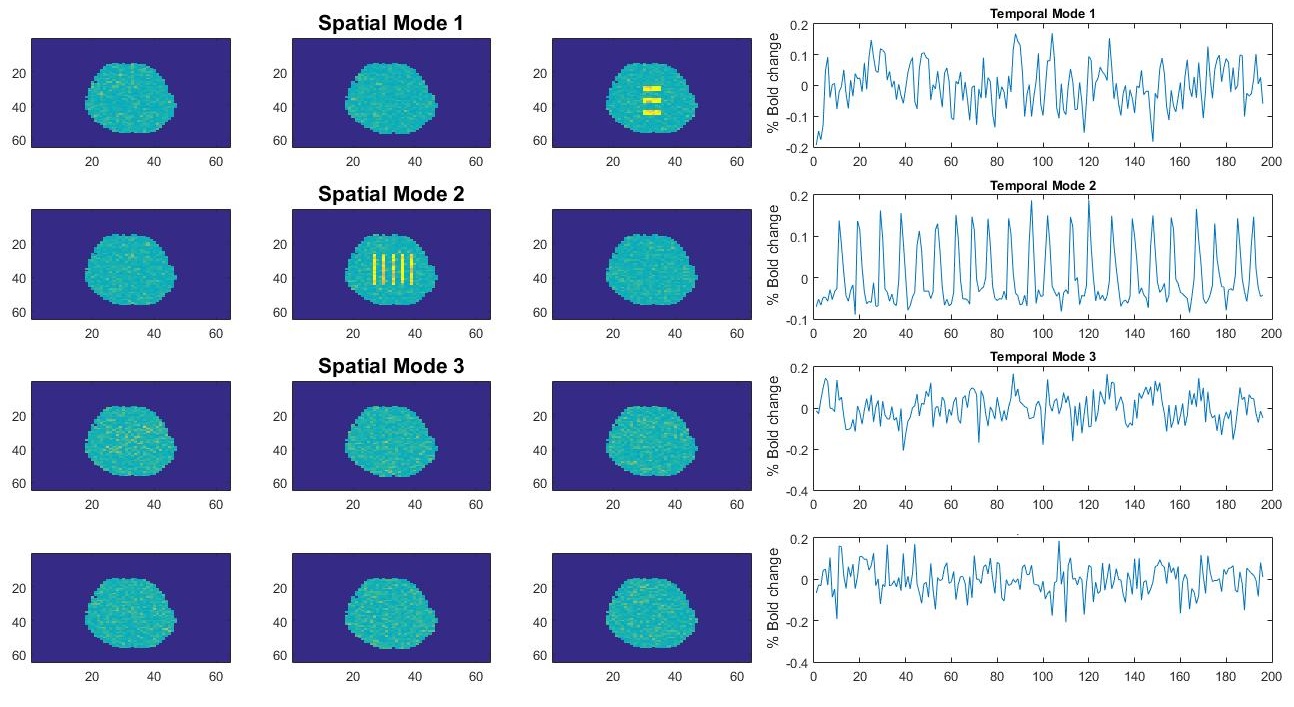}
  \caption [width=.80\linewidth]{TPICA of dataset G with SNR equal to 0.08.}
  \label{fig:tpicG8}
\end{figure}

\begin{figure}
  \centering
  \captionsetup{justification=centering}
  \includegraphics[width=1\linewidth]{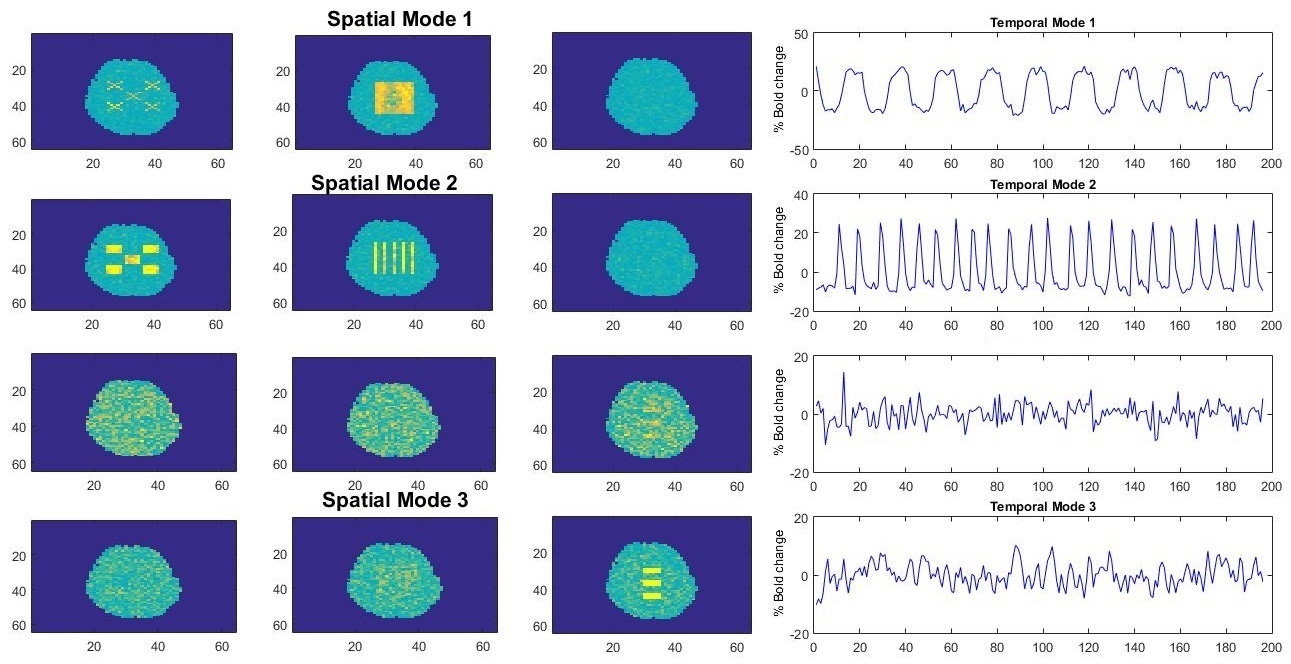}
  \caption [width=.80\linewidth]{CPD of dataset G with SNR equal to 0.15.}
  \label{fig:cpdG15}
\end{figure}

\begin{figure} 
  \centering
  \captionsetup{justification=centering}
  \includegraphics[width=1\linewidth]{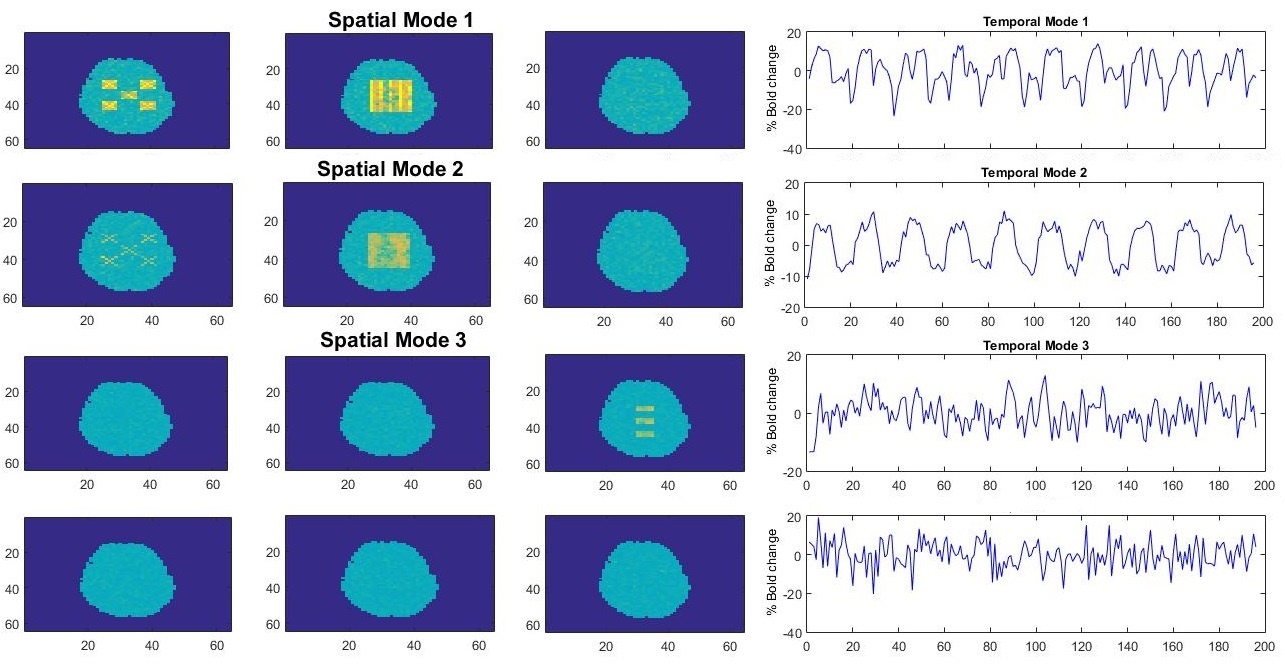}
  \caption [width=.80\linewidth]{CPD of dataset G with SNR equal to 0.08.}
\label{fig:cpdG8}
\end{figure}

\begin{figure}
  \centering
  \captionsetup{justification=centering}
  \includegraphics[width=1\linewidth]{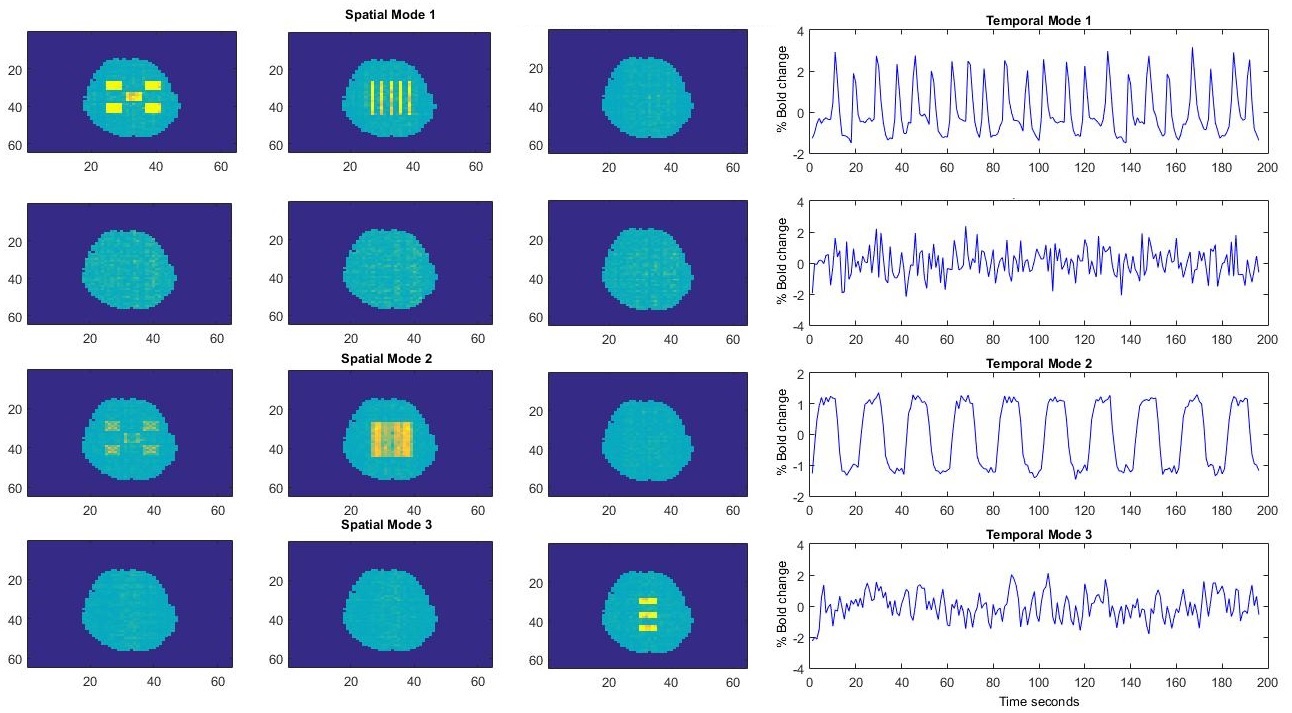}
  \caption [width=.80\linewidth]{BTD of dataset G with SNR equal to 0.15.}
  \label{fig:btdG15}
\end{figure}

\begin{figure}
  \centering
  \captionsetup{justification=centering}
  \includegraphics[width=.95\linewidth]{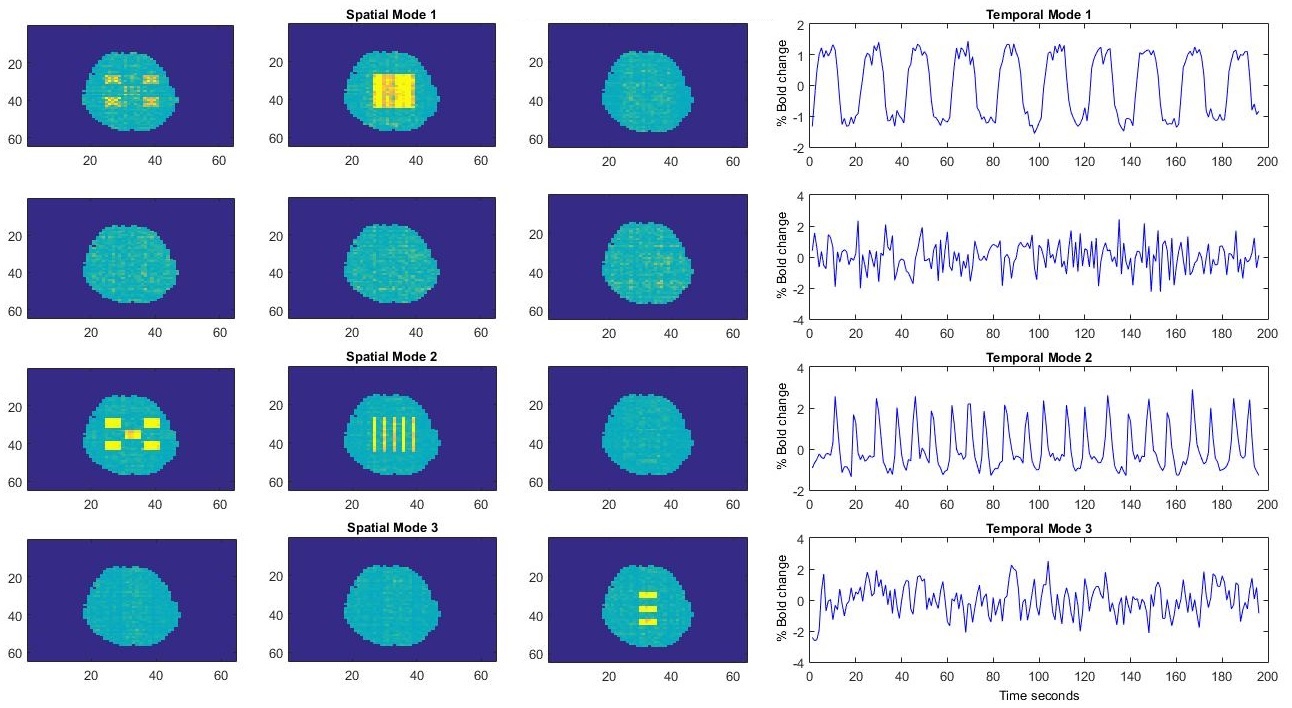}
  \caption [width=.80\linewidth]{BTD of dataset G with SNR equal to 0.08.}
\label{fig:btdG8}
\end{figure}

Figs.~\ref{fig:tpicG15}-\ref{fig:btdG8} validate some of the findings of the correlation tables. As mentioned before, even with SNR=0.15, TPICA fails to distinguish correctly the sources, as the second slice of the spatial maps of the first two sources is the same. CPD and BTD separate correctly the sources at SNR=0.15 (Figs. \ref{fig:cpdG15} and \ref{fig:btdG15}). At SNR=0.08, CPD can not distinguish the second slice of the first two sources (similarly to TPICA) (Fig.~\ref{fig:cpdG8}), while BTD still gives an almost perfect result (Fig.~\ref{fig:btdG8}). The rank of BTD, $L_r$, is set to 3 but, as mentioned in the first experiment, the method is insensitive to overestimation of the rank.

 \begin{table} [ht]
  \centering
  \begin{adjustbox}{width=1\linewidth}
  \begin{tabular}{||c||c|c|c||c|c|c||c|c|c||}
  \hline
   \multirow{2}{*}{Maps} & \multicolumn{3}{c}{SNR1=0.15} & \multicolumn{3}{|c|}{SNR2=0.12} & \multicolumn{3}{|c|}{SNR3=0.08*} \\
    \hhline{~---------}
    & Map 1 & Map 2 & Map 3 & Map 1 & Map 2 & Map 3 & Map 1 & Map 2 & Map 3  \\
  \hline
  \hline
   \centering TPICA M1 & 0.85 & 0.75 & 0.11 & 0.65 & 0.47 & 0.12 & 0.55 & 0.29 & 0.14\\ 
  TPICA M2& 0.69 & 0.89 & 0.09 & 0.59 & 0.48 & 0.12 & 0.54 & 0.32 & 0.11 \\ 
  TPICA M3 & 0.11& 0.12 & 0.98 & 0.12 & 0.12 & 0.94 & 0.20 & 0.19 & 0.88  \\
  \hline
     \centering CPD M1 & 0.96 & -0.20 & -0.12 & 0.87 & 0.44 & 0.21 & 0.75 & -0.47 & 0.22\\ 
  CPD M2 & -0.48 & 0.94 & 0.12 & 0.51 & 0.90 & -0.12 & -0.54 & 0.72 & 0.29  \\ 
  CPD M3 & -0.11 & 0.11 & 0.98 & 0.24 & 0.11 & 0.94 & 0.25 & 0.28 & 0.90 \\
  \hline  
  \centering BTD M1 & 0.96 & -0.15 & 0.1 & 0.94 & 0.23 & 0.12 & 0.88 & -0.32 & 0.13 \\ 
  BTD M2 & -0.25 & 0.94 & 0.18 & 0.28 & 0.92 & 0.18 & -0.27 & 0.88 & 0.19\\ 
  BTD M3 & -0.12 & 0.17 & 0.95 & 0.13 & 0.18 & 0.95 & 0.14 & 0.2 & 0.92 \\
  \hline
\end{tabular}
\end{adjustbox}
  \caption{Correlation of spatial maps of dataset G.}
  \label{tab:dg}
  \vspace{-2mm}
\end{table}

 \begin{table}
  \centering
  \begin{adjustbox}{width=1\linewidth}
  \begin{tabular}{||c||c|c|c||c|c|c||c|c|c||}
  \hline
   \multirow{2}{*}{Time Courses } & \multicolumn{3}{c}{SNR=0.15} & \multicolumn{3}{|c|}{SNR=0.12} & \multicolumn{3}{|c|}{SNR=0.08*} \\
    \hhline{~---------}
    & TC 1 & TC 2 & TC 3 & TC 1 & TC 2 & TC 3 & TC 1 & TC 2 & TC 3  \\
  \hline
  \hline
   \centering TPICA T1 & 0.96 & -0.54 & 0.13& 0.54 & 0.43 & 0.2 & 0.40 & 0.14 & 0.2\\ 
  TPICA T2 & -0.31 & 0.80 & 0.11 & 0.45 & 0.50 & -0.13 & 0.38 & 0.23 & -0.16  \\ 
  TPICA T3 & 0.14 & -0.12 & 0.94 & 0.20 & 0.24 & 0.93 & 0.22 & -0.50 & 0.60   \\
  \hline
   \centering CPD T1 & 0.96 & -0.24 & -0.11& 0.84 & 0.26 & 0.2 & 0.62 & -0.50 & -0.2\\ 
  CPD T2 & -0.30 & 0.92 & 0.21 & 0.32 & 0.840 & -0.23 & -0.38 & 0.53 & -0.26  \\ 
  CPD T3 & 0.11 & 0.11 & 0.94 & 0.14 & 0.14 & 0.93 & -0.22 & 0.3 & 0.65   \\
  \hline
   \centering BTD T1 & 0.96 & 0.21 & 0.11& 0.94 & 0.22 & 0.12 & 0.88 & 0.32 & 0.18\\ 
  BTD T2 & 0.20 & 0.92 & -0.11 & -0.22 & 0.890 & -0.13 & 0.23 & 0.83 & -0.16  \\ 
  BTD T3 & 0.1 & 0.11 & 0.96 & -0.14 & 0.15 & 0.91 & 0.21 & -0.29 & 0.85   \\
  \hline
  \multicolumn{10}{l}{*TPICA in SNR=0.08 often fails to separate 3 components and gives only two} 
\end{tabular}
\end{adjustbox}
  \vspace{-3mm}
  \caption{Correlation of time courses of dataset G.}
  \label{tab:tg}
\end{table}

\clearpage

\subsection{Simulation with different types of artifacts} 

\begin{figure} [ht]
\centering
\captionsetup{justification=centering}
\includegraphics[width=.95\textwidth]{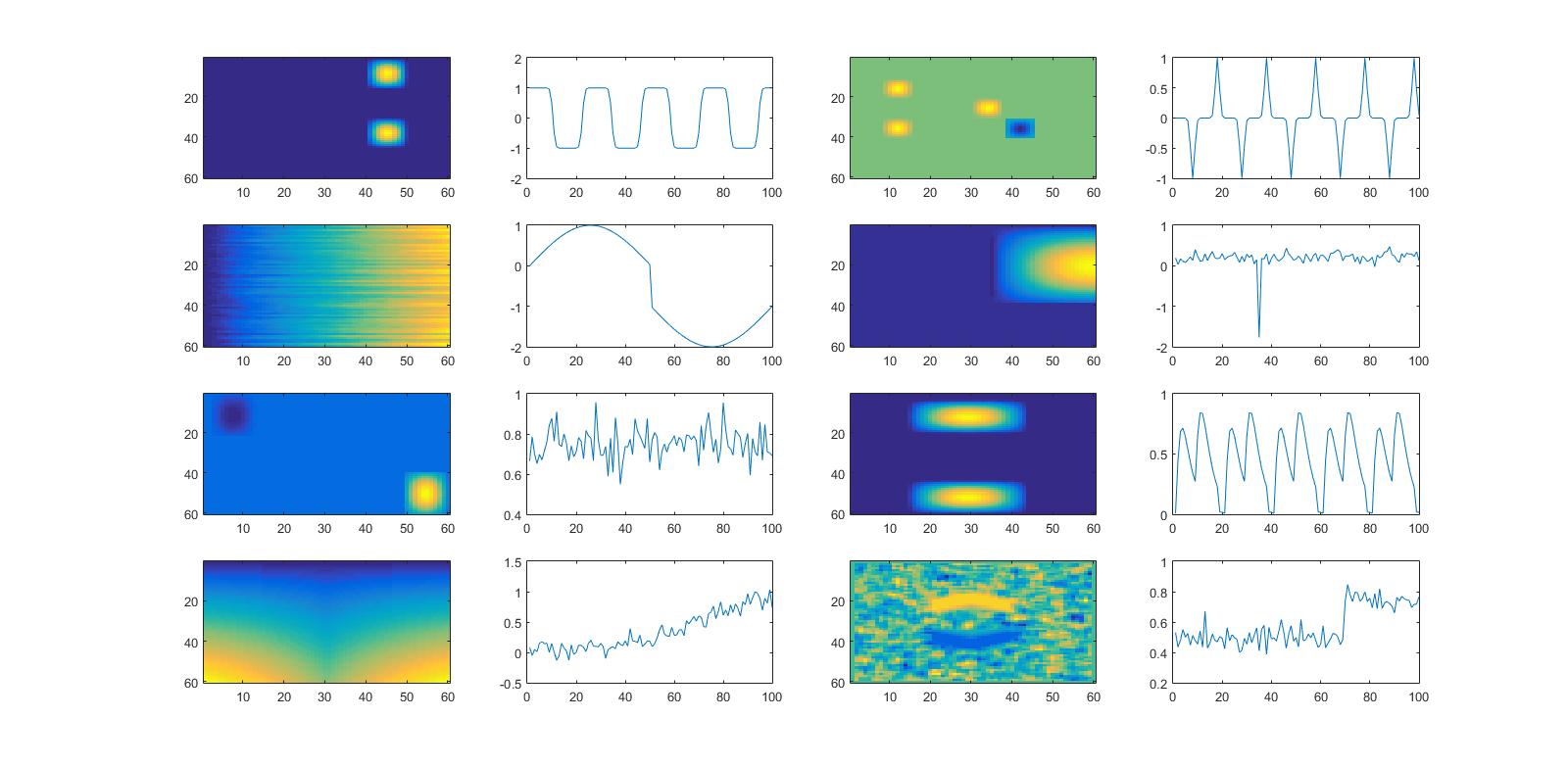}
\caption{\label{fig:data}  Data set from Erhardt et al. \cite{2011_erhardt_comparison}.}
\end{figure}

An fMRI-like dataset \cite{2011_erhardt_comparison} (Fig.~\ref{fig:data}) representing three sources of interest and five artifacts is used as a common set of sources for different subjects in this experiment. Among the eight spatial maps, one is task-related (1), two are transiently task-related (2, 6), and five are artifact-related (3, 4, 5, 7, 8). Furthermore, five of those are super-Gaussian (1, 2, 5, 6, 8), two are sub-Gaussian (3, 7), and one is Gaussian (4).

Besides the existing artifacts (3, 4, 5, 7, 8), and similarly to the previous two experiments, Gaussian noise was added for testing the three different methods. The experiment was carried out, first, without the Gaussian source (source 4, Fig.~\ref{fig:data}), which had the highest rank (equal to 52) and the highest energy; so it is the artifact which influences our result most.

\newpage
\subsubsection{Seven sources}

The components extracted in all decomposition methods in the first part of the experiment are seven while the rank $L_r$ of the BTD was selected equal to five.

\begin{figure}
\centering
\captionsetup{justification=centering}
\includegraphics[width=0.73\textwidth]{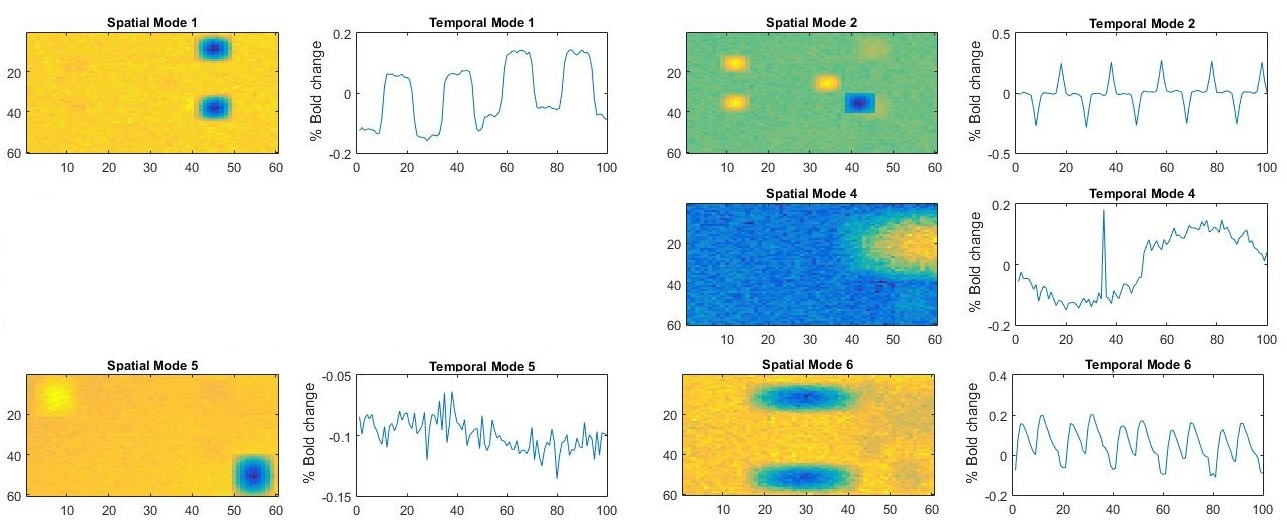}
\caption{\label{fig:tp1ad}  TPICA with SNR=0.8.}
\vspace{-2mm}
\end{figure}

\begin{figure} [ht]
\centering
\captionsetup{justification=centering}
\includegraphics[width=0.79\textwidth]{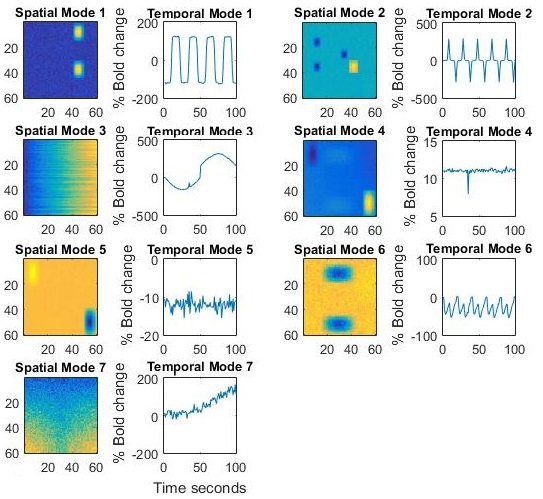}
\caption{\label{fig:cpd1ad}  CPD decomposition with SNR=0.8.}
\end{figure}

\begin{figure}[h]
\centering
\captionsetup{justification=centering}
\includegraphics[width=0.8\textwidth]{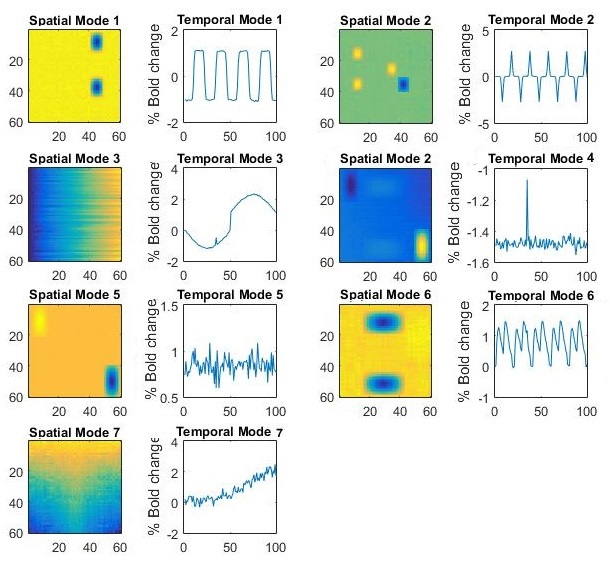}
\caption{\label{fig:btd1ad}  BTD decomposition with SNR=0.8.}
\end{figure}

\begin{figure} [h]
\centering
\captionsetup{justification=centering}
\includegraphics[width=0.8\textwidth]{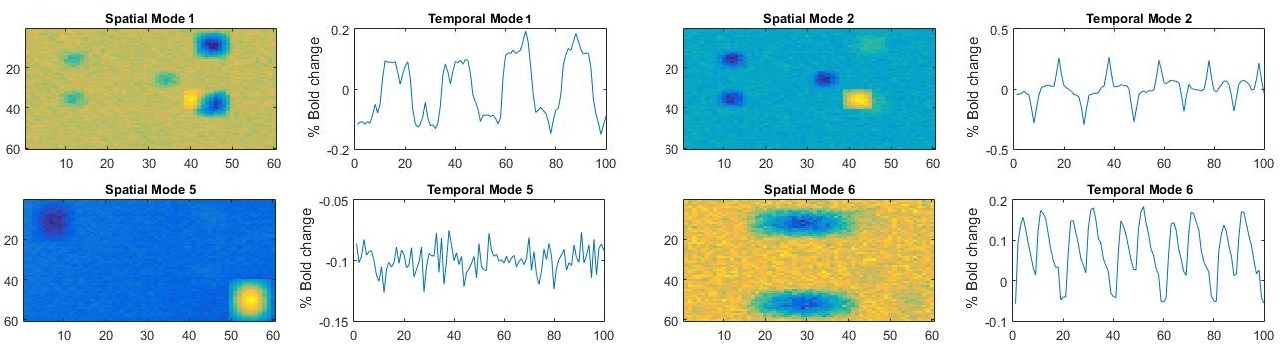}
\caption{\label{fig:tp2ad}  TPICA with SNR=0.4.}
\end{figure}

\begin{figure}
\centering
\captionsetup{justification=centering}
\includegraphics[width=0.7\textwidth]{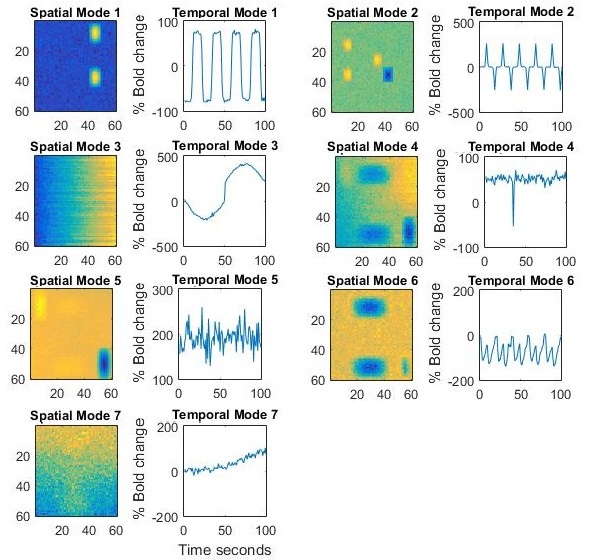}
\caption{\label{fig:cpd2ad} CPD decomposition with SNR=0.4.}
\vspace{-3mm}
\end{figure}

\begin{figure}
\centering
\captionsetup{justification=centering}
\includegraphics[width=0.7\textwidth]{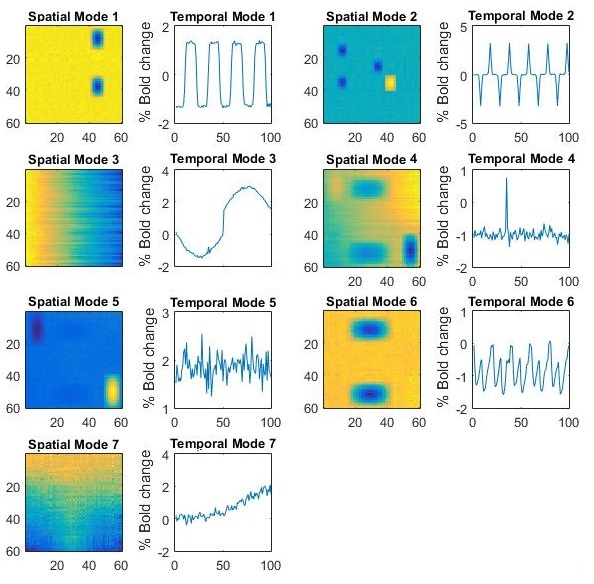}
\caption{\label{fig:btd2ad} BTD decomposition with SNR=0.4.}
\vspace{-3mm}
\end{figure}

Figs.~\ref{fig:tp1ad}-\ref{fig:btd2ad} exhibit the result of the decompositions, when TPICA is used as the method of the decomposition only five sources are distiniguished at SNR=0.8, while with SNR=0.4, only four out of the seven (no matter the number of the iterations of the algorithm). Furthermore, with higher noise (Fig.~\ref{fig:tp2ad}), even the task-related sources (sources~1, 2, 6), which are of the highest interest, are not well discriminated. In areas where sources~1 and 2 have a small spatial overlap, crosstalk between spatial maps can be noted.

\newpage

Both CPD and BTD at SNR=0.08 (Figs.~\ref{fig:cpd1ad}, \ref{fig:btd1ad}) result in an accurate decomposition. They separate all the sources (but one) almost perfectly. They both fail in distinguishing source 4 (which is the weakest in power). At SNR=0.04, the result of BTD remains stable and accurate (Fig.~\ref{fig:btd2ad}), while, on the other hand, the spatial maps, resulted from CPD, start becoming noisy and exhibiting some small crosstalk in the case of spatial map of source 6 (Fig.~\ref{fig:cpd2ad}).

\newpage

\clearpage
\subsubsection{Eight sources}

 When the last source (source~8) is added, the rank of the BTD needs to increase (e.g., $L_r>20$), since the source is of high rank and of high energy (or else the decomposition will fail similarly to Fig.~\ref{fig:examp}). This would be really time consuming and would increase the complexity of our decomposition. Instead of increasing the rank of all the matrices and have $R$ matrices $\textbf{X}_r \in \mathbb{R}^{I_1\times L_r}$ and $\textbf{Y}_r \in \mathbb{R}^{I_{2}I_{3}\times L_r}$, we can make a prediction about the amount of high energy-high rank sources (an overestimation does not affect our decomposition) and instead of BTD as in (\ref{btd4d}), we can have an adapted version of BTD:
 
 \begin{equation}
\label{eqadap}
\mathbf{\mathcal{T}}= \sum_{r=1}^{n} (\mathbf{X}_r {\mathbf{Y}_r^T}) \circ \mathbf{b}_r \circ \mathbf{c}_r + \sum_{r=n+1}^{R}  (\widetilde{\mathbf{X}}_r \widetilde{\mathbf{Y}}_r^{T}) \circ \mathbf{b}_r \circ \mathbf{c}_r,
\end{equation}

\begin{figure} 
\centering
\captionsetup{justification=centering}
\includegraphics[width=1\textwidth]{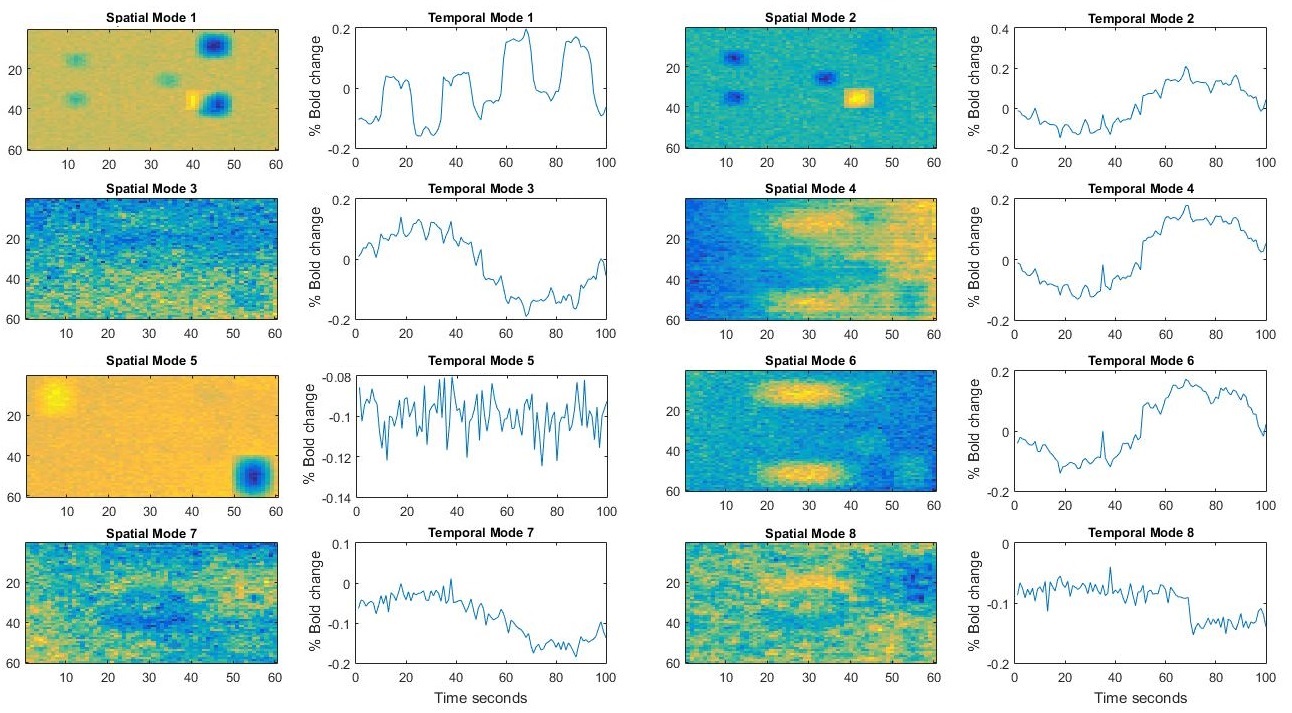}
\caption{\label{fig:tp28}  TPICA with SNR=0.4.}
\end{figure}
\noindent with $n$, $\textbf{X}_r$ and $\textbf{Y}_r$, matrices of rank $L_r$ and $R-n$, $\widetilde{\mathbf{X}}_r, \widetilde{\mathbf{Y}}_r$, matrices of rank $\widetilde{L}_r>L_r$.\footnote{The uniqueness of this adapted BTD is not studied here.}. For the needs of the experiments, 6 sources with $L_r=5$ (similarly to the seven-source experiment) were selected, and 2 with $\widetilde{L}_r=30$. A further increase to four sources with BTD rank equal to 30 will be also tested.

\begin{figure}
\centering
\captionsetup{justification=centering}
\includegraphics[width=0.95\textwidth]{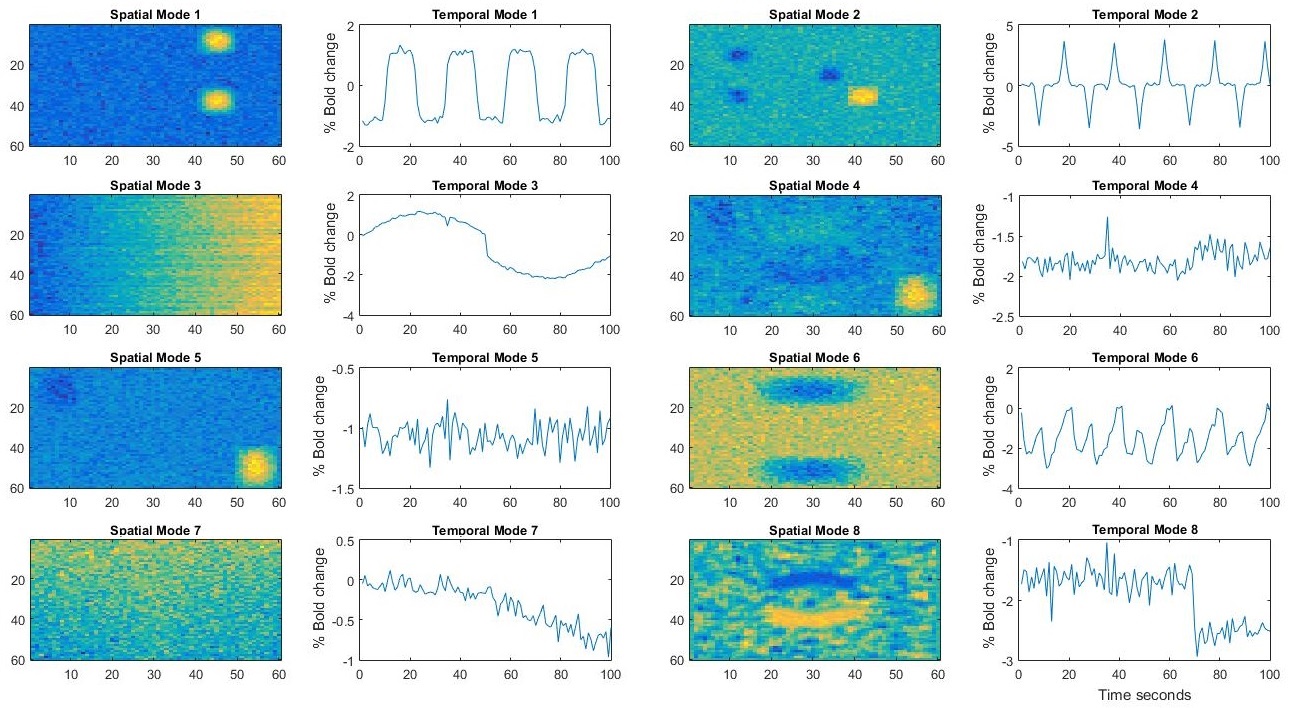}
\caption{\label{fig:cpd28}  CPD decomposition with SNR=0.4.}
\end{figure}

TPICA fails again to find correctly all the sources (Fig.~\ref{fig:tp28}),  only sources 1, 2, 5 and 8 can be correctly distinguished. Besides the fact that source 6, which is not an artifact, is not recognized, the first two (non-artifact) sources have crosstalk in their spatial maps (similarly to the case of 7 sources). CPD exhibits quite good results (Fig.~\ref{fig:cpd28}); it recognizes all the non-artifact sources (1, 2, 6) correctly, although with some noise, and it can also find sources 3, 4 and 8 and the time course of source 7. The spatial map of source 7 is not well recognized, while, once again the fourth source, which is weak, is not detected at all.

\begin{figure} 
\centering
\captionsetup{justification=centering}
\includegraphics[width=0.98\textwidth]{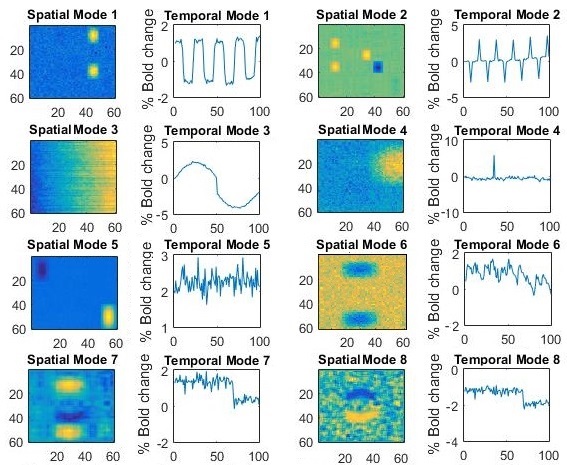}
\caption{\label{fig:bt28}  Adapted BTD decomposition with 2 sources having high rank.}
\end{figure}

\begin{figure} 
\centering
\captionsetup{justification=centering}
\includegraphics[width=0.82\textwidth]{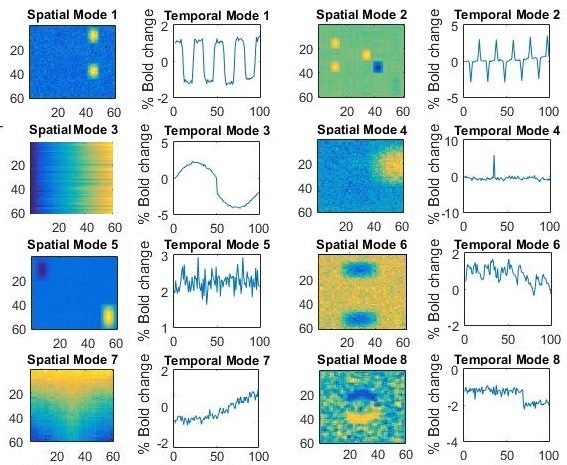}
\caption{\label{fig:bt28_4}  Adapted BTD decomposition with 4 sources having high rank.}
\end{figure}

The results of the adapted BTD (Figs.~\ref{fig:bt28} and \ref{fig:bt28_4}) are really interesting. First of all, in Fig.~\ref{fig:bt28}, for the first time source 4 is recognized correctly. The model ``decides'' to give the freedom of higher rank to source 8 (which is of the highest rank) and source 4 (which could not be recognized previously). If the number of sources with rank $L_r=30$ is increased to four, source 7, which was not previously recognized, is now correctly separated (Fig.~\ref{fig:bt28_4}). 

\newpage

\clearpage

\section{Tensors in Single-Subject Experiments}
Tensor decomposition methods have been introduced in the field of fMRI data analysis, for effectively performing blind source separation in multi-subject (or multi-trial) experiments. The dimensions of the tensor in most of these methods are voxels $\times$ time $\times$ subject/trial, and hence in single-subject, single-trial cases a multiway tensor can not be formed and matrix decomposition methods (ICA \cite{2008_xiong_ica}, k-SVD \cite{2015_zhao_supervised}, etc.) are used instead. In the proposed spatially folded BTD approach, the extra dimension of subject or trial is not necessary in order to give rise to a tensor, since it can be formed from the intrinsic spatial dimensions (voxels X $\times$ voxels YZ $\times$ time) (Fig.~\ref{fig:4d}). Results from a tensorial fMRI analysis of single-subject scans will be reported elsewhere.
\vspace{-2mm}

\section{ Conclusions}

In this paper, a novel approach, compared to the classical unfolding of the original 3D fMRI brain imaging, has been presented. The new path is based on a higher-dimensional unfolding, in an effort to exploit better the original 3D (tensor)- spatial structure of the data. This leads to 4-way tensors which makes it necessary  to employ the BDT method for tensor decomposition. Extensive simulation results demonstrated the enhanced robustness in the presence of noise of the new method compared to CPD-based decompositions.  In cases of spatial and temporal overlap, both CPD and BTD give better results than TPICA, albeit at the cost of a higher sensitivity w.r. to the rank estimate is observed (as reported in \cite{2007_stegeman_comparing} and \cite{2013_helwig_critique} also). 

In terms of computational resources, the complexity of BTD is higher than CPD (no  compression is assumed). It was also observed that when the rank, $L_r$, of the BTD increases, the result of the decomposition is equally good or even better, at the expense of higher complexity. In the rare case of an activation source, which has high rank and high energy simultaneously, BTD needs the rank, $L_r$, to be higher (and hence higher complexity) an enhanced source separation potential. The use of an adapted BTD procedure can help in preserving low computational cost and good decomposition results.

\section{Acknowledgment}

The authors would like to thank A. Stegeman and N. Helwig for providing the datasets used in~\cite{2007_stegeman_comparing} and~\cite{2013_helwig_critique}, respectively.

\bibliographystyle{myIEEEtran}
\bibliography{My_Library}

\end{document}